\documentclass[11pt]{article}
\usepackage{amssymb,amsmath,epsf}
\newtheorem{theorem}{Theorem}[section]

\newtheorem{theorem-definition}[theorem]{Theorem-Definition}
\newtheorem{theorem-construction}[theorem]{Theorem-Construction}
\newtheorem{lemma-definition}[theorem]{Lemma--Definition}
\newtheorem{proposition-definition}[theorem]{Proposition--Definition}
\newtheorem{lemma}[theorem]{Lemma}
\newtheorem{proposition}[theorem]{Proposition}
\newtheorem{corollary}[theorem]{Corollary}
\newtheorem{conjecture}[theorem]{Conjecture}
\newtheorem{definition}[theorem]{Definition}

\begin{document}
\newcommand{\Z}{{\mathbb Z}}
\newcommand{\R}{{\mathbb R}}
\newcommand{\Q}{{\mathbb Q}}
\newcommand{\C}{{\mathbb C}}
\newcommand{\lms}{\longmapsto}
\newcommand{\lra}{\longrightarrow}
\newcommand{\hra}{\hookrightarrow}
\newcommand{\ra}{\rightarrow}
\newcommand{\sgn}{\rm sgn}
\begin{titlepage}
\title{Cluster ${\cal X}$-varieties, amalgamation  
 and Poisson-Lie groups}  
\author{V.V. Fock, A. B. Goncharov}
\end{titlepage}
\date{}
\maketitle
\qquad \qquad \qquad \qquad \qquad \quad 
{\it To Vladimir Drinfeld  for his 50th birthday}

\tableofcontents

 \vskip 6mm \noindent

\section{Introduction} 

{\bf 1.1 Summary}. 
In this paper, starting from a split semisimple real Lie group $G$ with trivial center, we  define a family of varieties with additional structures. We describe them as the {\em cluster  ${\cal X}$-varieties}, as defined in \cite{FG2}. In particular they are Poisson varieties. We define canonical Poisson maps of these varieties to the group $G$ equipped with the standard Poisson-Lie structure defined by V.Drinfeld in (\cite{D}, \cite{D1}). One of them maps to the group birationally and thus provides $G$ with canonical rational coordinates. 

We introduce a simple but important operation of {\it amalgamation} of cluster varieties. 
 Our varieties are obtained as amalgamations of certain {\it elementary cluster varieties} ${\cal X}_{\mathbf J(\alpha)}$, assigned to positive simple roots  $\alpha$ of the root system of $G$. An elementary cluster variety ${\cal X}_{\mathbf J(\alpha)}$ is a split algebraic torus of dimension $r+1$, where $r$ is the rank of $G$. Its cluster, and in particular  Poisson, structure is described in a very simple way by the Cartan matrix for $G$.  Since one of them is a Zariski open part of $G$, we can  develop  the Poisson-Lie group structure on $G$ from scratch, without  the $r$-matrix formalism, getting as a benefit canonical (Darboux) coordinates for the Poisson structure on $G$. Some of our varieties are very closely related to the double Bruhat cells studied by A.Berenshtein, S.Fomin and A.Zelevinisky in \cite{FZ}, \cite{BFZ3}, \cite{BZq}. 

Using quantization of cluster ${\cal X}$-varieties developed in Section 4 of (\cite{FG2}) we get as a byproduct a quantization (i.e. a noncommutative $q$-deformation) of our varieties. The quantum version of the 
operation of amalgamation generalizes the standard quantum group structure. 

The results of this paper enter as a building block into a description of the cluster 
structure of the moduli spaces of local systems on surfaces studied in \cite{FG1}. 

\vskip 3mm
{\bf 1.2 Description of the results}. 
We start the paper with a brief recollection of the definition and properties of cluster ${\cal X}$-varieties. Let us briefly discuss some of their features, postponing the detailed discussion till Section 2. 

\vskip 3mm
\paragraph{1.2.1 Cluster ${\cal X}$-varieties.} 
 Cluster ${\cal X}$-varieties are determined by combinatorial data similarly (although different in some details) to the one used for definition of cluster algebras in \cite{FZI}, that is, by a  {\em cluster seed} ${\mathbf I}$, which is a quadruple $(I, I_0, \varepsilon, d)$, where

i) $I$ is a finite set;

ii)  $I_0 \subset I$ is its subset; 

iii) $\varepsilon$ is a matrix $\varepsilon_{ij}$, where $i,j \in I$, such that 
$\varepsilon_{ij} \in {\mathbb Z}$ unless $i,j \in I_0$.

iv) $d = \{d_i\}$, where $i \in I$, is a set of positive integers, such that the matrix $\widehat{\varepsilon}_{ij}=\varepsilon_{ij}d_j$ is skew-symmetric. 

The elements of the set $I$ are called {\em vertices}, the elements of $I_0$ are called {\em frozen vertices}.
 
Given a seed $\mathbf I$, every non-frozen vertex $k \in I-I_0$ gives rise to a {\it mutation}, producing a new, mutated seed $\mu_k(\mathbf I)$. Compositions of mutations are called {\it cluster transformations of seeds}.

\vskip 3mm 
Following Section 2 of \cite{FG2},  we associate to a seed ${\mathbf I}$ a torus ${\mathcal X}_{\mathbf I} = ({\mathbb G}_m)^I$ with a Poisson structure given by 
$$
\{x_i,x_j\}=\widehat{\varepsilon}_{ij}x_ix_j
$$ where $\{x_i| i\in I\}$ are the standard coordinates on the factors. We shall call it the {\em seed ${\cal X}$-torus}. 
Cluster transformations of seeds give rise to Poisson birational transformations between the seed tori, called cluster transformations. 
Gluing the seed ${\cal X}$-tori according to these birational transformations we get a scheme ${\cal X}_{|\mathbf I|}$ over $\Z$, called below {\it cluster ${\cal X}$-varieties}. (However ${\cal X}_{|\mathbf I|}$  may not be a scheme of finite type, and thus ${\cal X}_{|\mathbf I|}\otimes \Q$  may not be an algebraic variety.) 

\vskip 3mm 
In \cite{FZI}, the values $\varepsilon_{ij}$ for $i\in I_0$ were not defined, so $\varepsilon_{ij}$ was a rectangular matrix with integral entries. In our approach the frozen variables play an important role. The values $\varepsilon_{ij}$, when $i,j\in I_0$, are essential, and may not be integers. Let us elaborate this point. 

\vskip 3mm
\paragraph{1.2.2  Amalgamation.} We introduce operations of {\it amalgamation} and {\it defrosting} of seeds. The amalgamation of a collection of  seeds $\mathbf I(s)$, parametrised by a set $S$, is a new seed $\mathbf K = (K, K_0, \varepsilon_{ij}, d)$. 
The set $K$ is defined by gluing some of the frozen vertices of the sets $I(s)$. The frozen subset $K_0$ is obtained by gluing the frozen subsets $I_0(s)$. The rest of the data of $\mathbf K$ is also inherited from 
the ones of $\mathbf I(s)$. Defrosting simply shrinks the subset of the frozen vertices of $\mathbf K$, without changing the set $K$.  One can defrost any subset of $K$ such that $\varepsilon_{ij} \in \Z$ for any $i,j$ from this subset. All seeds in our paper are obtained by amalgamation followed by defrosting of certain elementary seeds. All vertices of the elementary seeds are frozen. 

\vskip 3mm
In order to state our results we have to introduce some notation related to the group. 
 
\vskip 3mm
\paragraph{1.2.3 Notation.} 
Let $G$ be a semisimple adjoint Lie group of rank $r$. There is the following data associated to $G$: the set of positive simple roots $\Pi$, the Cartan matrix $C_{\alpha\beta}=2\frac{(\alpha,\beta)}{(\alpha,\alpha)}$, where $\alpha,\beta \in \Pi$, and the multipliers $d^\alpha=(\alpha,\alpha)/2, \alpha \in \Pi$, such that the matrix $\widehat{C}_{\alpha\beta} = C_{\alpha\beta}d^\beta$ is symmetric. 
Let $\Pi^-$ be the set of negative simple roots and let ${\mathfrak W}$ be the semigroup freely generated by $\Pi$ and $\Pi^-$. Any element $D$ of ${\mathfrak W}$ is thus a word $\mu_1\cdots\mu_{l(D)}$ of the letters from the alphabet $\Pi\cup\Pi^-$, where $l(D)$ is its length. For $\alpha\in \Pi$ we shall denote by $\bar{\alpha}$ the opposite element from $\Pi^-$.

\vskip 3mm
\paragraph{1.2.4 The braid and Hecke semigroups.} 
Let ${\mathfrak B}$ be the quotient of the semigroup ${\mathfrak W}$ by
\begin{equation}\label{c1}
\alpha\bar{\beta}=\bar{\beta}\alpha
\end{equation}
\begin{equation}
\begin{array}{lclcl}\label{c2}
\alpha\beta\alpha = \beta\alpha\beta &\mbox{and}& \bar{\alpha}\bar{\beta}\bar{\alpha}=\bar{\beta}\bar{\alpha}\bar{\beta} &\mbox{if}& C_{\alpha\beta}=C_{\beta\alpha}=-1,\\
\alpha\beta\alpha\beta = \beta\alpha\beta\alpha &\mbox{and}& \bar{\alpha}\bar{\beta}\bar{\alpha}\bar{\beta}=\bar{\beta}\bar{\alpha}\bar{\beta}\bar{\alpha} &\mbox{if}& C_{\alpha\beta}=2C_{\beta\alpha}=-2,\\
\alpha\beta\alpha\beta\alpha\beta = \beta\alpha\beta\alpha\beta\alpha &\mbox{and}& \bar{\alpha}\bar{\beta}\bar{\alpha}\bar{\beta}\bar{\alpha}\bar{\beta}=\bar{\beta}\bar{\alpha}\bar{\beta}\bar{\alpha}\bar{\beta}\bar{\alpha} &\mbox{if}& C_{\alpha\beta}=3C_{\beta\alpha}=-3.
\end{array}
\end{equation}
The semigroup ${\mathfrak B}$ is called the {\it braid semigroup}. We denote by $p:{\mathfrak
W}\rightarrow{\mathfrak B}$ the canonical projection.

Another semigroup appropriate in our context is the further quotient of ${\mathfrak B}$ by the relations
\begin{equation}
\alpha^2 = \alpha; \quad \bar{\alpha}^2=\bar{\alpha},
\end{equation} 
denoted by ${\mathfrak H}$ and called the {\it Hecke semigroup}. It is isomorphic as a set to the square of the Weyl group of $G$. We call an element of ${\mathfrak W}$ or ${\mathfrak B}$ {\em reduced} if its length is minimal among the elements having the same image in ${\mathfrak H}$.

\vskip 3mm 
Now we are ready to discuss our main goals and results. 

\vskip 3mm
\paragraph{1.2.5 Cluster ${\cal X}$-varieties related to the braid semigroup and their properties.} 
We will define a cluster ${\cal X}$-variety ${\cal X}_B$ associated to
 any element $B\in {\mathfrak B}$. For this purpose, given a $D \in {\mathfrak W}$,  we will define a seed 
$\mathbf J(D)=(J(D),J_0(D),\varepsilon(D),d(D))$. We prove that the seeds corresponding to different elements of $p^{-1}(B)$ are related by cluster transformations. Moreover, we define the evaluation and multiplication maps 
$$
ev: {\cal X}_B \rightarrow G, \qquad 
m:{\cal
X}_{B_1}\times{\cal X}_{B_2}\rightarrow {\cal X}_{B_1B_2}$$
The correspondence $B \mapsto {\cal X}_B$ and the maps $m$ and $ev$ are to satisfy the following properties:
\begin{enumerate}
\item $ev$ is a Poisson map.
\item $m$ is an  amalgamation followed by defrosting of cluster ${\cal X}$-varieties (and thus a Poisson map).
\item The multiplication maps are associative in the obvious sense. 
\item Multiplication commutes with the evaluation, i.e. 
the following diagram is commutative:
\begin{equation} \label{mevcomm1}
\begin{array}{ccc}
{\cal X}_{B_1}\times {\cal X}_{B_2}&\stackrel{(ev, ev)}{\lra}& G \times G\\
\downarrow m &&\downarrow m\\
{\cal X}_{B_1B_2}&\stackrel{ev}{\lra}&G
\end{array}
\end{equation}
\end{enumerate}
We would like to stress that the multiplication $m$ is a projection with fibers of nonzero dimension. 

\vskip 3mm
\paragraph{1.2.6 Cluster ${\cal X}$-varieties related to the Hecke semigroup.}
Let $\pi: {\mathfrak B} \to {\mathfrak H}$ be the canonical projection of semigroups. Considered as a projection of sets it has  a canonical  splitting 
$s: {\mathfrak H} \to {\mathfrak B}$. 
Namely, for every $H \in {\mathfrak H}$ 
there is a unique reduced element $s(H)$ in $\pi^{-1}(H)$,
 the {\it reduced representative} of $H$ 
in ${\mathfrak B}$. So given an element 
$H \in {\mathfrak H}$ there is  a cluster variety ${\cal X}_{s(H)}$. 
Abusing notation, we will denote it by ${\cal X}_{H}$. 

A rational map of cluster ${\cal X}$-varieties  is a {\it cluster projection} if 
in a certain cluster coordinate system it is obtained by forgetting one or more cluster coordinates. 

We show that

\begin{enumerate}
\item There is a canonical cluster projection $\pi: {\cal X}_{B} \to {\cal X}_{\pi(B)}$. 
By the very definition, it  is an isomorphism if $B$ is reduced. 
\item There is a multiplication map 
$
m_{\cal H}:{\cal
X}_{H_1}\times{\cal X}_{H_2}\rightarrow {\cal X}_{H_1H_2}
$, defined as the composition
$$
{\cal X}_{H_1}\times{\cal X}_{H_2} := {\cal X}_{s(H_1)}\times{\cal X}_{s(H_2)} \stackrel{m}{\lra}
{\cal X}_{s(H_1)s(H_2)} \stackrel{\pi}{\lra} {\cal X}_{H_1H_2}
$$
So it  is a composition of an amalgamation, defrosting, and cluster projection.  
It follows from (\ref{mevcomm1}) that  the maps $m_{\cal H}$ and  $m$ are related by a commutative diagram
\begin{equation} \label{acons1}
\begin{array}{ccc}
{\cal X}_{B_1}\times{\cal X}_{B_2} &\stackrel{m}\lra &{\cal X}_{B_1B_2}\\
\downarrow \pi \times \pi &&\downarrow \pi \\
{\cal X}_{H_1}\times{\cal X}_{H_2} &\stackrel{m_{\cal H}}\lra &{\cal X}_{H_1H_2}
\end{array}
\end{equation}
where the vertical maps are the canonical cluster projections. 
\item The multiplication maps are associative in the obvious sense 
(see the Remark below).

\item If $H \in \mathfrak H$, the map $ev: {\cal X}_{H} \hra G$ is injective 
at the generic point. If $H$ is the longest element of ${\mathfrak H}$, then the image of $ev$ is Zariski dense in $G$. The map 
$ev: {\cal X}_{B} \to G$  is a composition 
$$
{\cal X}_{B} \stackrel{\pi}{\lra} {\cal X}_{H} \stackrel{ev_{\cal H}}{\hra} G, \qquad H = \pi (B)
$$
\end{enumerate}

\vskip 3mm
{\bf Remark 1}. A part of the above data are axiomatized as follows. 
Given a semigroup $\mathfrak S$, we assign to every $s \in \mathfrak S$ 
an object ${\cal X}_s$ of a monoidal category ${\cal M}$ (e.g., 
the category 
of Poisson varieties with the product as a monoidal structure), and for every 
pair $s, t \in \mathfrak S$ 
a canonical morphism
 $m_{s,t}: {\cal X}_s \times {\cal X}_t \lra {\cal X}_{st}$. They must satisfy 
an associativity constraint, i.e. for every $r,s,t \in \mathfrak S$, the following diagram is commutative:  
\begin{equation} \label{acons}
\begin{array}{ccc}
{\cal X}_{r}\times{\cal X}_{s}\times{\cal X}_{t} &
\stackrel{{\rm Id} \times m_{s,t}}{\lra} &{\cal X}_{r}\times{\cal X}_{st} \\
m_{r,s}\downarrow \times {\rm Id}&& \downarrow m_{r, st}\\
{\cal X}_{rs}\times{\cal X}_{t} & \stackrel{m_{rs, t}}{\lra} & {\cal X}_{rst} 
\end{array}
\end{equation}

\vskip 3mm
{\bf Remark 2}. 
Recall Lusztig's coordinates on the group $G$ (\cite{L1}).
 Let $E^{\alpha}(t)$ and $F^{\alpha}(t)$ be the two standard 
one parametric subgroups 
corresponding to a simple root $\alpha$. Denote by 
$X^{\alpha}(t)$ the element $E^{\alpha}(t)$ if $\alpha \in \Pi$ and 
$F^{\alpha}(t)$ if $\alpha \in \Pi^-$. 
A reduced decomposition of the longest element in $W \times W$ is encoded by 
 a sequence 
$\alpha_1, ..., \alpha_{2m}$ of 
$2m$ elements of $\Pi \cup \Pi^-$, where $2m+r = {\rm dim}G$. 
There is a birational isomorphism 
\begin{equation} \label{acons2}
H \times {\mathbb G}_m^{2m}\lra G, \quad (H, t_1, ..., t_{2m})\lms 
H X^{\alpha_1}(t_1), ..., X^{\alpha_{2m}}(t_{2m}) 
\end{equation}
There are similar coordinates on all double Bruhat cells (\cite{FZ}). 
However they are not cluster ${\cal X}$-coordinates. 
Our coordinates are related to them by monomial 
transformations.

\vskip 3mm
\paragraph{1.2.7 Quantization.} Cluster ${\cal X}$-varieties were quantized in 
Section 4 of \cite {FG2}.  
The operations of amalgamation, defrosting and cluster projection have 
straightforward generalizations to the quantum ${\cal X}$-varieties. 
Thus we immediately get $q$-deformations 
of the considered above cluster ${\cal X}$-varieties for the braid and  
Hecke semigroups. 

We understood the category of quantum spaces as in loc. cit. 
So in particular 
a morphism of quantum cluster spaces ${\cal X}^q \to {\cal Y}^q$, by definition, 
is given by a compatible collection of morphisms of the 
corresponding quantum tori algebras going in the opposite direction, i.e. 
the  
${\cal Y}$-algebras map to the corresponding ${\cal X}$-algebras. 

It follows that the quantum spaces enjoy 
the properties similar to the  listed above properties 
of their classical counterparts:  

\begin{enumerate}
\item There is a canonical projection $\pi: {\cal X}_{B} \to {\cal X}_{\pi(B)}$.

\item  There are multiplication 
maps of quantum spaces, for the braid and Hecke semigroups:  
$$
m^q: {\cal
X}^q_{B_1}\times{\cal X}^q_{B_2}\rightarrow {\cal X}^q_{B_1B_2}, \qquad 
m^q_{\cal H}: {\cal
X}^q_{H_1}\times{\cal X}^q_{H_2}\rightarrow {\cal X}^q_{H_1H_2},  
$$

\item  They are related by the $q$-version of the diagram (\ref{acons1}), and 
satisfy the associativity constraints 
given by the $q$-versions of the diagram 
(\ref{acons}). 

\end{enumerate}

{\bf Remark}. One can show that there is a quantum evaluation map to the 
quantum deformation (\cite{D1}) of the algebra of regular functions 
of $G$. Unlike the other properties, this is not 
completely straightforward, and will be elaborated elsewhere. 

\vskip 3mm
\paragraph{1.2.8 Proofs.} They are easy if  
the Dynkin diagram of $G$ is simply-laced. 
The other cases are reduced to the rank two cases.  
The  $B_2$ case 
can be done by obscure calculations.
The $G_2$ case is considerably more difficult.  

A more conceptual approach is provided by the operation of {\it cluster folding} \cite{FG2}, 
 briefly reviewed in Subsection 3.6, which clarifies the picture in the $B_2$ case and seems to be indispensable in the $G_2$ case.

\vskip 3mm
\paragraph{1.2.9 Cluster structures of moduli spaces 
of triples of flags 
of types $A_3$ and $G_2$.} In the process of proof  we work with 
cluster ${\cal X}$-varieties related to the moduli spaces ${\rm Conf}_3({\cal B}_{B_2})$  and ${\rm Conf}_3({\cal B}_{G_2})$ of configurations of 
triples of flags in the Lie groups of type $B_2$ and $G_2$, respectively. There is a canonical embedding 
${\rm Conf}_3({\cal B}_{B_2}) \hra {\rm Conf}_3({\cal B}_{A_3})$ provided by the folding of the latter.  
It allows to reduce the study of the former to the study of the latter. 
In the two appendices we investigate the cluster structures of the moduli spaces 
${\rm Conf}_3({\cal B}_{A_3})$ and ${\rm Conf}_3({\cal B}_{G_2})$ in detail.  
Here is what  we learned. 

Recall the moduli 
space ${\cal M}_{0, 6}$ of configurations  of $6$ points on ${\Bbb P}^1$. It has an 
${\cal X}$-cluster structure of finite type $A_3$.  In the Appendix 2 we 
 construct 
an  explicit birational isomorphism 
$$
\Phi: {\rm Conf}_3({\cal B}_{A_3}) \stackrel{\sim}{\lra} {\cal M}_{0, 6}
$$
respecting the ${\cal X}$-cluster structures.

\vskip 3mm
The investigation of the moduli space ${\rm Conf}_3({\cal B}_{G_2})$ 
turned out to be a subject of independent interest, which reveals the following 
story, discussed in the Appendix 1.

We say that two seeds $\mathbf I = (I, I_0, \varepsilon, d)$ and 
$\mathbf I' = (I', I'_0, \varepsilon', d')$ are isomorphic, if there is a set isomorphism $\varphi: I \to I'$, preserving the frozen vertices and the functions $\varepsilon$ and $d$-functions.  

\begin{definition} \label{8.20.05.1}
A cluster ${\cal X}$-variety is of {\em $\varepsilon$-finite} type if the set of the isomorphism classes of its seeds is finite. 
\end{definition}

Any cluster ${\cal X}$-variety gives rise to an orbifold, called the {\it modular orbifold}, see Section 2 in \cite{FG2}. We recall its definition in Section 4.1. Its dimension is the dimension of the cluster ${\cal X}$-variety minus one. 
In the $\varepsilon$-finite case the modular orbifold is glued from a finite number of simplices. It is non-compact, unless the cluster ${\cal X}$-variety is of finite cluster type. Here is the main result: 

\begin{theorem} \label{3.9.05.1a}
a) The cluster ${\cal X}$-variety corresponding to the moduli space ${\rm Conf}_3({\cal B}_{G_2})$ is of $\varepsilon$-finite type. 
The number of  the isomorphism classes of its  seeds is seven. 

b) The corresponding modular orbifold is a manifold. It is homeomorphic to $S^3 - L$, where $L$ is a two-component link, and $\pi_1(S^3-L)$ is isomorphic to the braid group of type $G_2$.  
\end{theorem}

It is well known that the complement to the discriminant variety of type $G_2$ in $\C^2$ is a $K(\pi, 1)$-space where $\pi$ is the braid  group of type $G_2$. Its intersection with a sphere $S^3$ containing the origin has two connected components. We conjecture that it is isomorphic to $S^3 - L$. 

The {\it mapping class group} of a cluster ${\cal X}$-variety  was defined in loc. cit. It acts by automorphisms of the cluster ${\cal X}$-variety. 
It is always infinite if the cluster structure is of 
$\varepsilon$-finite, but not of finite type. Theorem \ref{3.9.05.1a} immediately 
implies 
the following:

\begin{corollary}\label{3.9.05.1b}
The mapping class group of the cluster ${\cal X}$-variety corresponding to 
${\rm Conf}_3({\cal B}_{G_2})$ 
is an infinite quotient of 
the braid group of type $G_2$.
\end{corollary}
Conjecturally it coincides with the braid group. 
This is the first example of an infinite mapping class group 
different from the mapping class groups of surfaces. 

\vskip 3mm 
 {\bf Acknowledgments}. V.F. was supported by the grants CRDF 2622; 2660. He is grateful to the 
Math departments of Universite Lyon-1 and Brown University for
hospitality, and to Lacrimiora Iancu for a stimulating lecture. 
A.G. was supported by the NSF grant DMS-0400449. He is grateful to IHES for the hospitality. 

\section{Cluster ${\cal X}$-varieties and amalgamation.}

In this Section we recall some definitions from \cite{FG2}. 

\vskip 3mm
{\bf 2.1 Basic definitions}. A {\em cluster seed}, or just {\em seed}, ${\mathbf I}$ is a quadruple $(I, I_0, \varepsilon, d)$, where

i) $I$ is a finite set;

ii)  $I_0 \subset I$ is its subset; 

iii) $\varepsilon$ is a matrix $\varepsilon_{ij}$, where $i,j \in I$, such that $\varepsilon_{ij} \in {\mathbb Z}$ unless $i,j \in I_0$.

iv) $d = \{d_i\}$, where $i \in I$, is a set of positive integers, such that the matrix $\widehat{\varepsilon}_{ij}=\varepsilon_{ij}d_j$ is skew-symmetric. 

The elements of the set $I$ are called {\em vertices}, 
the elements of $I_0$ are called {\em frozen vertices}.  
The matrix $\varepsilon$ is called a {\em cluster function}, 
the numbers $\{d_i\}$ are called {\em multipliers}, and the function 
$d$ on $I$ whose value at $i$ is $d_i$ is called a {\it multiplier function}. 
We omit $\{d_i\}$ if all of them are equal to one, and therefore the matrix 
$\varepsilon$ is skew-symmetric, and we omit the set $J_0$ if it is empty. 

The seed ${\mathbf I}=(I, I_0, \varepsilon, d)$ is called a {\em subseed} of the seed ${\mathbf I}'=
(I',I'_0,\varepsilon, d)$ if $I\subset I'$, $I_0\subset I'_0$ and the functions $\varepsilon$ and $d$ are the restriction of $\varepsilon'$ and $d'$, respectively. In this case we denote ${\mathbf I}$ by ${\mathbf I}'\vert_I$.

Recall the multiplicative group scheme ${\mathbb G}_m$. 
It is defined as the spectrum of the ring $\Z[X,X^{-1}]$. 
The direct product of several copies of the multiplicative group is called a split algebraic torus, or simply a torus. 
The readers who are not used to the language of schemes may just fix, once forever, a field $K$, 
and replace everywhere ${\mathbb G}_m$ by $K^\times$; indeed ${\mathbb G}_m(K)= K^\times$. 

For a seed ${\mathbf I}$ we associate a torus ${\mathcal X}_{\mathbf I} = ({\mathbb G}_m)^I$ 
with a Poisson structure given by 
\begin{equation}
\{x_i,x_j\}=\widehat{\varepsilon}_{ij}x_ix_j \label{Poisson_bracket}
\end{equation} where $\{x_i| i\in I\}$ are the 
standard coordinates on the factors. We shall call it the {\em seed ${\cal X}$-torus}. 

Let ${\mathbf I}=(I, I_0,\varepsilon, d)$ and ${\mathbf I}'=(I', I'_0,\varepsilon', d')$  be two seeds, and $k\in
I$. A {\em mutation in the vertex $k$} is an isomorphism $\mu_k: I\rightarrow I'$ satisfying the following conditions:
\begin{enumerate}
\item $\mu_k(I_0)=I'_0$,
\item $d'_{\mu_k(i)}=d_i$,
\item $\varepsilon'_{\mu_k(i)\mu_k(j)}=
\left\{ \begin{array}{lll} -\varepsilon_{ij} & \mbox{ if } & i=k \mbox{ or } j=k\\
\varepsilon_{ij} + \varepsilon_{ik}\max(0,\varepsilon_{kj}) \mbox{ if } \varepsilon_{ik}\geq 0\\
\varepsilon_{ij} + \varepsilon_{ik}\max(0,-\varepsilon_{kj}) \mbox{ if } \varepsilon_{ik}<0
\end{array}\right.$
\end{enumerate}

A {\em symmetry} of a seed ${\mathbf I}=(I, I_0,\varepsilon, d)$ is an  automorphism $\sigma$ of the set $I$ preserving the subset $I_0$, the matrix $\varepsilon$ and the numbers $d_i$. In other words it satisfies the conditions:
\begin{enumerate}
\item $\sigma(I_0)=I_0$,
\item $d_{\sigma(i)}=d_i$,
\item $\varepsilon_{\sigma(i) \sigma(j)}=\varepsilon_{ij}$
\end{enumerate}

Symmetries and mutations induce (rational) maps between the corresponding  
seed ${\cal X}$-tori, which are  denoted by the same symbols $\mu_k$ and $\sigma$ and given by the formulae

$$x_{\sigma(i)}=x_i$$
and
$$
x_{\mu_k(i)} = \left\{\begin{array}{lll} x_k^{-1}& \mbox{ if } & i=k \\
    x_i(1+x_k)^{\varepsilon_{ik}} & \mbox{ if } & \varepsilon_{ik}\geq 0 \mbox{ and } i\neq k \\
 x_i(1+(x_k)^{-1})^{\varepsilon_{ik}} & \mbox{ if } & \varepsilon_{ik}\leq 0 \mbox{ and } i\neq k
\end{array} \right..
$$

A {\em cluster transformation} between two seeds (and between two seed
${\cal X}$-tori) is a composition of symmetries and mutations. If the source and the target of a cluster transformation coincide we call this map a {\em cluster automorphism}. Two seeds are called {\em equivalent} if they are related by a cluster transformation. The equivalence class of a seed ${\mathbf I}$ is denoted by $|{\mathbf I}|$

Thus we have defined two categories. The first one have cluster seeds as objects and cluster transformation as morphisms. The second one have seed ${\mathcal X}$-tori as objects and cluster transformations of them 
as morphisms. There is a canonical functor from the first to the second. 
The objects in these two categories are the same.
However there are more  morphisms in the second category. 

A {\em cluster ${\cal X}$-variety} is obtained by taking a union of all seed ${\cal X}$-tori related to a given seed ${\mathbf I}$ by cluster transformations, and gluing them together  using the above birational isomorphisms. It is denoted by ${\cal X}^{}_{|{\mathbf I}|}$. Observe that
 the cluster ${\cal X}$-varieties corresponding to equivalent seeds are isomorphic.  Every particular seed ${\cal X}$-torus  provide our cluster variety with a rational  coordinate system. The corresponding rational functions are called {\em cluster coordinates}.

Since in the sequel we shall extencively use compositions of mutations we would like to introduce a 
shorthand notation for them. Namely, we denote an expression $\mu_{\mu_{i}(j)}\mu_{i}$ by $\mu_j\mu_k$, $\mu_{\mu_{\mu_{i}(j)}\mu_{i}(k)}\mu_{\mu_{i}(j)}\mu_{i}$ by $\mu_k\mu_j\mu_i$, and so on. We will also say that two sequences of mutations are equivalent ($\stackrel{\sim}{=}$) if they coincide as maps between the $\mathcal X$-tori up to permutation of coordinates. 

The cluster transformations have the following basic properties, see \cite{FG2}, Section 2:
\begin{enumerate}
\item Every seed is related to other seeds by exactly $\sharp (I- I_0)$ mutations.
\item Cluster transformations form a groupoid. In particular the inverse of a mutation is a mutation: $\mu_k\mu_k=id$. Cluster automorphisms form a group called {\em mapping class group}. The groups of cluster automorphisms of equivalent seeds are isomorphic.
\item Cluster transformations preserve the Poisson structure. In  particular a cluster ${\cal X}$-manifold has a canonical Poisson structure and the automorphism group of this manifold acts on it by Poisson transformations.
\item Cluster transformations are given by rational functions with positive integral  coefficients.
\item If $\varepsilon_{ij}=\varepsilon_{ji}=0$ then $\mu_i\mu_j\mu_j\mu_i=id$.
\item If $\varepsilon_{ij}=-\varepsilon_{ji}=-1$ then $\mu_i\mu_j\mu_i\mu_j\mu_i\stackrel{\sim}{=}id$. (This is called the {\em pentagon relation}.)
\item If If $\varepsilon_{ij}=-2\varepsilon_{ji}=-2$ then $\mu_i\mu_j\mu_i\mu_j\mu_i\mu_j\stackrel{}{=}id$.
\item If If $\varepsilon_{ij}=-3\varepsilon_{ji}=-3$ then $\mu_i\mu_j\mu_i\mu_j\mu_i\mu_j\mu_i\stackrel{}{=}id$.
\end{enumerate}

Conjecturally all relations between mutations are exhausted by the properties 5-8.

\vskip 3mm
\paragraph{2.2 Amalgamation.} We start from the simplest example: 
the amalgamation of two seeds. Let $\mathbf J = (J,J_0,\varepsilon,d)$ and $\mathbf I=(I,I_0,\zeta,c)$ be two seeds and let $L$ be a set embedded into both $I_0$ and $J_0$ in a such a way that for any $i,j \in L$ we have 
 $c(i)=d(i)$. Then the amalgamation of  $\mathbf J$ and $\mathbf I$ is a seed $\mathbf K = (K,K_0,\zeta,b)$, such that $K=I \cup_L J$, $K_0 = I_0\cup_L J_0$ and
$$\zeta_{ij}=\left\{\begin{array}{lcl}
               0 &\mbox{if}& i\in I- L \mbox{ and } j\in J- L\\
               0 &\mbox{if}& i\in J- L \mbox{ and } j\in I- L\\
               \eta_{ij}&\mbox{if}&i\in I- L \mbox{ or } j\in I- L\\
               \varepsilon_{ij}&\mbox{if}&i\in J- L \mbox{ or } j\in J- L\\
               \eta_{ij}+\varepsilon_{ij}&\mbox{if}&i,j\in L
                    \end{array}
             \right.
$$
This operation induces a homomorphism ${\cal X}_{\mathbf J} \times {\cal X}_{\mathbf I} \to {\cal X}_{\mathbf K}$ between the corresponding seed ${\cal X}$-tori given by the rule
\begin{equation}\label{restrx}
z_i = \left\{\begin{array}{lcl}x_i&\mbox{if}&i\in I- L\\
                               y_i&\mbox{if}&i\in J-  L\\
                               x_iy_i&\mbox{if}&i\in L
              \end{array}
       \right.
\end{equation}

It is easy to check that it respects the Poisson structure and
 commutes with cluster transformations, 
thus is defined for the cluster ${\cal X}$-varieties, and not only for the seeds.

If there is a subset $L' \subset L$, such that 
$\varepsilon_{ij}+\eta_{ij}\in \mathbb Z$ when $i,j \in L$, then we can 
{\it defrost} the vertices of  $L'$, getting a new seed 
$(K,K_0- L',\zeta,b)$. This way we get a different cluster ensemble, since we can now mutate the elements of $L'$ as well.

\vskip 3mm
Now let us present the general definition. Let 
$$
\mathbf I({s})=
(I({s}), I_0({s}), \varepsilon({s}), d({s})), \qquad {s} \in S
$$
be a family of seeds parametrised by a set $S$. Let us 
glue the sets $I^{s}$ in such a  way  that 

\vskip 3mm
(a) only frozen vertices can be glued. 

(b) If $i \in I({s})$ and $j \in I({t})$ 
are glued, then $d({s})_{i} = d({t})_j$. 
\vskip 3mm

Let us denote by $K$ the set obtained by gluing the sets $I({s})$. 
\vskip 3mm
Alternatively, the {\it gluing data} can be described by the following data:
\vskip 3mm
(i) A set $K$.

(ii) A collection of injective maps $p_{s}: I({s}) \hra K$, 
$s \in S$, whose images cover $K$, 
two images may intersect only at the frozen elements. 

(iii) The multiplier function on $\cup_{s \in S}I({s})$ descends to a function $d$ on the 
set $K$. 
\vskip 3mm
In other words, there is  a cover of the set $K$ by the subsets $I({s})$, 
 any two elements covering the same point are frozen, 
and the values of the multiplier functions at these elements coincide. 
\vskip 3mm

We identify $\varepsilon({s})$ with a function on 
the square of the image of the set $I({s})$, and denote by 
$\varepsilon({s})'$ its extension by zero to $K^2$. 
Then we set 
\begin{equation} \label{restrxs}
\varepsilon:= \sum \varepsilon({s})'
\end{equation}
There is a map 
$
P:\cup_{s \in S}I({s}) \to K
$.  
We set
 $$
K_0:= P(\cup_{s \in S}I_0({s}))
$$
There is a unique 
function $c$ on $K$ such that 
$p_{s}^*d = d({s})$ 
for any $s \in C$.

\begin{definition} \label{} 
The seed $\mathbf K:= (K, K_0, \varepsilon, d)$ is  
the amalgamation of the seeds  
$\mathbf I({s})$ with respect to the given gluing data (i)-(iii). 
\end{definition}

\begin{lemma} \label{9.21.04.10}
The amalgamation of seeds commutes with cluster transformations.
 \end{lemma}

{\bf Proof}. Thanks to (ii), for any element 
$i \in I({s}) - I_0({s})$ 
one has $|P^{-1}(p_{s}(i))| =1$. Thus when we do a mutation 
in the direction $p_{s}(i)$, we can 
change the values of the cluster function only on the subset  
$p_{s}({I({s})}^2)$. The lemma follows. 

\vskip 3mm

\vskip 3mm
{\it The amalgamation map of cluster ${\cal X}$-varieties}. 
Let us consider the following  map 
of the ${\cal X}$-tori:  
\begin{equation} \label{2.14.05.10}
m: \prod_{s \in S}{\cal X}_{\mathbf I({s})} \to 
{\cal X}_{\mathbf K}, \qquad 
m^*x_i = \prod_{j\in P^{-1}(i)}x_j
\end{equation}

The following lemma is obvious:

\begin{lemma} \label{9.20.04.11} The maps (\ref{2.14.05.10})  commute with 
mutations, and thus give rise to a map of cluster ${\cal X}$-varieties, 
called the amalgamation map:
$$
m: \prod_{s \in S}{\cal X}_{|\mathbf I({s})|} \to 
{\cal X}_{|\mathbf K|}
$$  
\end{lemma} 
 
\vskip 3mm
\paragraph{Defrosting.} 
Let 
$L \subset K_0$. Assume 
that the function $\varepsilon$ restricted to 
$L\times L - K_0\times K_0$ 
takes values in $\Z$. Set 
$K'_0:= K_0  -  L$. 
Then there is a new seed  $\mathbf K':= (K, K_0', \varepsilon, c)$. We say that 
{\it the seed   $\mathbf K'$ 
is obtained from $\mathbf K$ by defrosting of  $L$}. 
There is a canonical open embedding
${\cal X}_{|\mathbf K|} \hra {\cal X}_{|\mathbf K'|}$. 

Amalgamation followed by defrosting is the key operation which we use below. 
Abusing notation, sometimes one may refer to this operation simply by amalgamation. 
However then the defrosted subset must be specified. 

\section{Cluster ${\cal X}$-varieties related to a group $G$}

\vskip 3mm
\paragraph{3.1 An example: rational coordinates on $PGL_2$.} Observe that one has 
 $$
H(x) := 
\begin{pmatrix}x^{1/2}&0\\0&x^{-1/2}\end{pmatrix} \ \stackrel{PGL(2,\mathbb C)}{=}\ 
\begin{pmatrix}x&0\\0&1\end{pmatrix} 
$$ 
They are elements of a Cartan subgroup in $PGL(2,\Bbb C)$. 
 
Consider a map $ev_{\bar \alpha \alpha}: ({\mathbb C}^\times)^3\rightarrow PGL(2,{\mathbb C})$ given by
$$
ev_{\bar \alpha \alpha}: (x_0,x_1,x_2) \mapsto
H(x_0)
\begin{pmatrix}1&0\\1&1\end{pmatrix}
H(x_1)
\begin{pmatrix}1&1\\0&1\end{pmatrix}
H(x_2) = 
(x_0x_1x_2)^{-1/2}\begin{pmatrix}x_0x_2(1+x_1)&x_0\\x_2& 1\end{pmatrix} 
$$

Consider also another map $ev_{\alpha \bar \alpha}: 
({\mathbb C}^\times)^3\rightarrow PGL(2,{\mathbb C})$:
$$
ev_{\alpha \bar \alpha}: (y_0,y_1,y_2) \mapsto
H(y_0) 
\begin{pmatrix}1&1\\0&1\end{pmatrix}
H(y_1)
\begin{pmatrix}1&0\\1&1\end{pmatrix}
H(y_2) = 
(y_0y_1y_2)^{-1/2}\begin{pmatrix}y_0y_1y_2&y_0y_1\\
                y_1y_2& 1 + y_1
\end{pmatrix}
$$

One can see that these maps enjoy the following properties:
\begin{enumerate}
\item $ev_{\bar \alpha, \alpha}(x_0,x_1,x_2)=ev_{\alpha, \bar \alpha}(x_0(1+x_1^{-1})^{-1},x_1^{-1}, 
x_2(1+x_1^{-1})^{-1})$.
\item $ev_{\alpha, \bar \alpha}(y_0,y_1,y_2)=ev_{\bar \alpha, \alpha}(y_0(1+y_1),y_1^{-1},y_2(1+y_1))$.
\item Both maps are open embeddings.
\item The standard Poisson-Lie 
bracket on the group $PGL(2,{\mathbb C})$ reads as 
$$\{x_0,x_1\}=x_0x_1;\qquad  \{x_1,x_2\}=x_1x_2;\ \quad\{x_0,x_2\}=0;$$
$$\{y_0,y_1\}=-y_0y_1; \quad\{y_1,y_2\}=-y_1y_2;\quad\{y_0,y_1\}=0;$$
\end{enumerate}

Therefore $ev_{\bar \alpha \alpha}$ and $ev_{\alpha \bar \alpha}$ provide the group variety $PGL(2,{\mathbb C})$ with two rational coordinate systems. The transition between these coordinates is given by a mutation and thus the union of the images of $ev_{\bar \alpha, \alpha}$ and $ev_{\alpha, \bar \alpha}$ is a cluster variety corresponding to two equivalent seeds $(\{0,1,2\},\{0,2\},\varepsilon)$ and $(\{0,1,2\},\{0,2\},\eta)$, where $\varepsilon = -\eta =\begin{pmatrix}
0&1&0\\-1&0&1\\0&-1&0\end{pmatrix}$ (columns correspond to the first index).

One can also consider the maps 
$$
ev_{\alpha}, ev_{\bar \alpha}: 
({\mathbb C}^\times)^2 \rightarrow PGL(2,{\mathbb C})  \mbox{ and } 
ev_{\emptyset}:{\mathbb C}^\times \rightarrow PGL(2,{\mathbb C})
$$
 given by:
$$
ev_{\alpha}: (z_0,z_1) \mapsto
H(z_0) 
\begin{pmatrix}1&1\\0&1\end{pmatrix}  
H(z_1) = (z_0z_1)^{-1/2}
\begin{pmatrix}z_0z_1& z_0\\0&1\end{pmatrix}
$$
$$
ev_{\bar \alpha}: (w_0,w_1) \mapsto
H(w_0) 
\begin{pmatrix}1&0\\1&1\end{pmatrix}
H(w_1) = (w_0w_1)^{-1/2}
\begin{pmatrix}w_0w_1&0
\\w_1&1\end{pmatrix}
$$
$$
ev_{\emptyset}: (t) \mapsto 
 t^{-1/2} \begin{pmatrix}t&0\\0&1\end{pmatrix}
$$

and satisfying the following properties:
\begin{enumerate}
\addtocounter{enumi}{5}
\item The images of $ev_{\alpha}$, $ev_{\bar \alpha}$, $ev_{\emptyset}$ 
and the union of the images of $ev_{\bar \alpha \alpha}$ and 
$ev_{\alpha \bar \alpha}$ are pairwise disjoint. 
The complement to their union in the whole group is of codimension two. It consists of the antidiagonal matrices. 
\item The images of $ev_{\alpha}$,$ev_{\bar \alpha}$,$ev_{\emptyset}$ 
are Poisson subvarieties with respect to the standard Drinfel-Jimbo poisoon-Lie structure on $PGL(2,\mathbb C)$. The Poisson bracket is given by:
$$
\{z_0,z_1\}=z_0z_1;\qquad  \{w_0,w_1\}=w_0w_1;
$$
\end{enumerate}

In the above constructions one can replace $\mathbb C^{\times}$ by ${\mathbb G}_m$, $PGL(2,\mathbb C)$ by the group scheme $PGL(2)$, and upgrade all maps to the maps of the corresponding schemes. 
\vskip 3mm

Our aim now is a generalization of this picture in two directions.  
We consider split semisimple adjoint groups of higher ranks, 
and  construct Poisson varieties which map to the group 
respecting the Poisson structure, but may not inject into the group.

\vskip 3mm 
Below we will give two alternative definitions of the seed $\mathbf J(D)$. The first one 
is computation free: we define first the elementary 
seeds $\mathbf J(\alpha)$ corresponding to simple roots 
$\alpha$, and then define $\mathbf J(D)$ as an amalgamated product of the elementary 
 seeds 
$\mathbf J(\alpha_1)$, ...,$\mathbf J(\alpha_n)$, where $D= \alpha_1 ... \alpha_n$, followed by defrosting of some of the frozen variables. It is presented in the subsection 2, and it is the definition which we  
use proving the main properties of our varieties in the subsection 4. The second definition is given by defining directly 
all components of the seed; its most important part is an
 explicit formula for the function $\varepsilon_{ij}$. The second definition 
is given in the subsection 3.

\vskip 3mm
\paragraph{3.2 The seed $\mathbf J(D)$.} Let us assume first that $D = \alpha$ is a simple positive root. 
Then we set
$$
J(\alpha)=  J_0(\alpha):= (\Pi - \{\alpha\}) \cup \{\alpha_-\} \cup \{\alpha_+\}
$$
where $\alpha'$ and $\alpha''$ are certain new elements. There is a {\it decoration} map
$$
\pi: J(\alpha)\lra \Pi, 
$$
which sends  $\alpha'$ and $\alpha''$ to $\alpha$, and is the identity map on $\Pi - \alpha$.

The collection of multipliers $\{d^{\alpha}\}$ gives rise to a function ${\cal D}$ on the set $\Pi$: ${\cal D}(\alpha):= d^{\alpha}$. 
We define the multipliers for $\mathbf J(\alpha)$ as the function ${\cal D} \circ \pi$ on $J(\alpha)$. 

Finally, the function $\varepsilon(\alpha)$ is defined as follows. Its entry $\varepsilon(\alpha)_{\beta \gamma}$ 
is zero unless one of the indexes is decorated by $\alpha$. 
Further, 
\begin{equation}
\varepsilon(\alpha)_{\alpha' \beta}= \frac{C_{\alpha \beta}}{2}, \quad 
\varepsilon(\alpha)_{\alpha'' \beta}= -\frac{C_{\alpha \beta}}{2}, \quad 
\varepsilon(\alpha)_{\alpha' \alpha''}= -1 \label{Elementary_epsilon}
\end{equation}

If $D = \bar \alpha$ is a negative simple root, we have a similar  set
$$
J(\bar \alpha)= J_0(\bar \alpha):= (\Pi^- - \{\bar \alpha\}) \cup 
\{\bar \alpha'\} \cup \{\bar \alpha''\},
$$
a similar decoration $\pi: J(\bar \alpha)\lra \Pi$, and similar multipliers ${\cal D} \circ \pi$ on $J(\bar \alpha)$. 
The cluster function is obtained by reversing the signs, using the obvious identification of $J(\bar \alpha)$ and $J(\alpha)$: 
$\varepsilon(\bar \alpha):= -\varepsilon(\alpha)$:
$$
\varepsilon(\bar \alpha)_{\bar \alpha' \bar \beta}= -\frac{C_{\alpha \beta}}{2}, \quad 
\varepsilon(\bar \alpha)_{\bar \alpha'' \bar \beta}= \frac{C_{\alpha \beta}}{2}, \quad 
\varepsilon(\bar \alpha)_{\bar \alpha' \bar \alpha''}= 1
$$  

The torus ${\cal X}_{\alpha}$ can also be defined as follows. 
Let $U_{\alpha}$ be the one parametric unipotent subgroup 
corresponding to the root $\alpha$. Then 
${\cal X}_{\alpha} = H\times (U_{\alpha}-\{0\}) = HU_{\alpha} -H$.

{\it The general case}. Observe that the subset of the elements decorated 
by a simple 
root has one element unless this root is $\alpha$, when there are two elements. 
There is a natural linear order on the subset of elements of $\mathbf 
J(\alpha)$ decorated by a given simple positive root $\gamma$: 
it is given by $(\alpha_-, \alpha_+)$ in the only nontrivial case when 
$\alpha = \gamma$. 
So for a given simple positive root $\gamma$ there are the {\it minimal} 
and the {\it maximal}
 elements decorated by $\gamma$. 

\begin{definition} \label{2.13.05.1} Let $D = \alpha_1 ... 
\alpha_n \in {\mathfrak W}$. Then  $J(D)$ is obtained 
by gluing the sets $J(\alpha_1), ..., J(\alpha_n)$ as follows. 
For every $\gamma\in \Pi$, 
and for every $i = 1, ..., n-1$,  we glue the maximal $\gamma$-decorated 
element of $J(\alpha_i)$ and the minimal 
$\gamma$-decorated element of $J(\alpha_{i+1})$. 

The seed ${\widetilde {\mathbf J}}(D)$ is the amalgamated product for this  
gluing data of the 
seeds $\mathbf J(\alpha_1)$, ..., $\mathbf J(\alpha_n)$. 
\end{definition}

The seed $\widetilde {\mathbf J}(D)$ has the frozen vertices only: 
${\widetilde {J}}(D) = {\widetilde {J_0}}(D)$. To define the seed $\mathbf J(D)$ 
we will defrost some of them, making 
the set $J_0(D)$ smaller.

In Definition \ref{2.13.05.1} we glue 
only the elements decorated by the same positive simple root. Thus 
the obtained set $J(D)$ has a natural decoration $\pi: J(D) \lra \Pi$, extending the ones of the subsets $J(\alpha_i)$. 
Moreover, for every $\gamma\in \Pi$, the subset of $\gamma$-decorated elements of $J(D)$ has a natural linear order, 
induced by the ones on $J(\alpha_i)$, and the linear order of the word $D$. 

\begin{definition} \label{2.13.05.2}
The subset $J_0(D)$ is the union, over $\gamma \in \Pi$, of the extremal (i.e. minimal and maximal) elements for the defined above linear order on the $\gamma$-decorated part of $J(D)$. 

The seed ${{\mathbf J}}(D)$ is obtained from the seed 
${\widetilde {\mathbf J}}(D)$ by reducing $\widetilde J_0(D)$ to the subset 
$J_0(D)$. 
\end{definition}
Observe that $\varepsilon_{\alpha \beta}$ 
is integral unless both $\alpha$ and $\beta$ are in $J_0(D)$. 
Thus the integrality condition for $\varepsilon_{\alpha \beta}$ holds. 

\vskip 3mm
\paragraph{3.3 An alternative definition of the seed $\mathbf J(D)$.}
{\it The sets $J_0(D)$ and $J(D)$.} 
Given a positive simple root $\alpha\in \Pi$, denote by $n^\alpha(D)$ the number of occurrences of $\alpha$ and 
$\bar{\alpha}$ in the word $D$. We set 
$$
J^\alpha(D):= \{(^\alpha_i)|\alpha \in \Pi, 0\leq i \leq n^\alpha\}, \qquad 
J^\alpha_0(D):= \{(^\alpha_0)\}\cup \{(^\alpha_{n_\alpha})\}
$$
Then $J(D)$ (resp. $J_0(D)$) is the disjoint union of $J^\alpha(D)$ (resp. $J^\alpha_0(D)$) for all $\alpha \in \Pi$. 
Observe  that if a root $\alpha$ does not enter the word $D$ then $J^\alpha_0(D)=J^\alpha(D)$ is a one element set.

One can picture elements of the set $J(D)$ as associated to the  intervals between walls made by $\alpha$, $\bar{\alpha}$ or the ends of the word $D$ for some root $\alpha$. If at least one wall is just the end $D$, the corresponding element of $J(D)$ belongs to $J_0(D)$. We shall denote these elements by braces connecting the walls with the name of the corresponding elements in the middle.

\vskip 3mm
  {\bf Example.} Let $r=\mathrm{rk}\ G=3$ and $\Pi=\{\alpha,\beta,\gamma\}$. Take
$D=\alpha\bar{\beta}\bar{\alpha}\bar{\alpha}\beta$. Then
$$n^\alpha(D)=3,\  n^\beta(D)=2,\ n^\gamma(D)=0$$
$$J(D) = \{(^\alpha_0),
(^\alpha_1),(^\alpha_2),(^\alpha_3),(^\beta_0),(^\beta_1),(^\beta_2),(^\gamma_0)\}, $$
$$J_0(D) = \{(^\alpha_0), (^\alpha_3),(^\beta_0),(^\beta_2),(^\gamma_0)\}, $$
In brace notations the set $J(D)$ can be shown as: 
$$ \underbrace{
\lefteqn{\overbrace{\phantom{\hat{\beta}}
\hspace{9ex}}^{(_0^\beta)}\overbrace{\phantom{\hat{\beta}}\hspace{16.7ex}}^{(_1^\beta)}\overbrace{\phantom{\hat{\beta}}}^{(_2^\beta)}}
\lefteqn{\underbrace{\phantom{y}}_{(_0^\alpha)}
\underbrace{\phantom{y}\hspace{10.5ex}}_{(_1^\alpha)}\underbrace{\phantom{y}\hspace{4.5ex}}_{(_2^\alpha)}
\underbrace{\phantom{y}\hspace{9.5ex}}_{(_3^\alpha)}}
\;\;\quad\alpha\qquad\bar{\beta}\qquad\bar{\alpha}\qquad\bar{\alpha}
\qquad\beta\quad\;\;\;}_{(_0^\gamma)}
$$

{\it  A description  of $\varepsilon$ and $d$.} 
In order to give an explicit formula for the matrix $\varepsilon_{(^\alpha_i)(^\beta_j)}$ we 
introduce more notation. Let $n^\alpha(k)$ be the number of letters $\alpha$ or $\bar{\alpha}$ among the first $k$ letters of the word $D$. Let $\mu_k$ be the $k$-th letter of $D$ and let $sgn(k)$ be $+1$ if $\mu_k\in \Pi$ and $-1$ otherwise. Let finally $|\mu_k|= sgn(\mu_k)\mu_k\in \Pi$.

\begin{definition} \label{9.23.04.1} Let $D$ be a word of ${\mathfrak W}$, Then: 

The multipliers are given by the rule $d_{(^\alpha_i)}(D)=d^\alpha$. 

The integers $\widehat \varepsilon_{(^\alpha_i)(^\beta_j)}$ are defined by the following formula
$$
\sum_{(^\alpha_i)(^\beta_j)}\widehat \varepsilon_{(^\alpha_i)(^\beta_j)}\frac{\partial}{\partial
\log x^\alpha_i}\wedge\frac{\partial}{\partial \log x^\beta_j} =
$$
$$
=\frac{1}{2}\sum_{k=1}^{l(D)}\sum_\alpha sgn(\mu_k)\widehat{C}_{\mu(k)\alpha}\frac{\partial}{\partial \log
x^\alpha_{n^\alpha(k)}}\wedge\left(\frac{\partial}{\partial \log x^{|\mu_k|}_{n^{|\mu_k|}(k)-1}}
-\frac{\partial}{\partial \log x^{|\mu_k|}_{n^{|\mu_k|}(k)}}\right).
$$
\end{definition}
\vskip 3mm

{\bf Remark}. One can check that in the case when $D$ is reduced,  
our function $\varepsilon_{ij}$ is related to the cluster 
function $b_{ij}$ defined in \cite{BFZ3} 
for the corresponding double Bruhat cell as follows. Recall that 
$b_{ij}$ is not defined if both $i,j$ are frozen variables. 
Other then that, the  
values of $b_{ij}$  turn out to be the same as for $\varepsilon_{ij}$.

\vskip 3mm 
It is easy to prove the following properties of the matrix $\varepsilon$:
\begin{trivlist}
\item[\ 1.] $\varepsilon_{(^\alpha_i)(^\beta_j)}$ is integral unless both $(^\alpha_i)$ and $(^\beta_j)$ are in $I_0$.
\item[\ 2.] For a given $(^\alpha_i) \in I$ the number of $(^\beta_j)\in I$, such that $\varepsilon_{(^\alpha_i)(^\beta_j)}\neq 0$ is no more than twice the number of $\beta \in \Pi$ such that $C_{\alpha\beta}\neq 0$. In particular this number never exceeds 8.
\item[\ 3.] The value of $\varepsilon_{(^\alpha_i)(^\beta_j)}$ is determined by the patterns of the walls in $D$, corresponding to $(^\alpha_i)$ and $(^\beta_j)$. The list of all possibilities is too large to give it explicitly, but we give just some of them --- the patterns to the left and the corresponding values of 
$\varepsilon_{(^\alpha_i)(^\beta_j)}$ to the right (stars mean any roots or word ends, compatible with the pattern):
\begin{eqnarray*}
\lefteqn{\hspace{1ex}\overbrace{\phantom{\beta}\hspace{5ex}}^{(^\alpha_i)}\overbrace{\phantom{\beta}\hspace{4ex}}^{(^\alpha_{i+1})}}
*\qquad\alpha\qquad*&\quad& 1\\
\lefteqn{\hspace{1ex}\overbrace{\phantom{\beta}\hspace{10ex}}^{(^\alpha_i)}}
\lefteqn{\hspace{7ex}\underbrace{\qquad \phantom{y}\hspace{6ex}}_{(^\beta_j)}}
*\qquad\alpha\qquad\beta\qquad*
&\quad& C_{\alpha\beta}\\
\lefteqn{\hspace{1ex}\overbrace{\phantom{\beta}\hspace{10ex}}^{(^\alpha_i)}}
\lefteqn{\hspace{7ex}\underbrace{\qquad \phantom{y}}_{(^\beta_j)}}
\alpha\qquad\beta\qquad\qquad \mbox{ or }\qquad
\lefteqn{\overbrace{\phantom{\beta}\hspace{3.5ex}}^{(^\alpha_0)}}
\lefteqn{\underbrace{\qquad \phantom{y}\hspace{5.5ex}}_{(^\beta_0)}}
\qquad\alpha\qquad\beta
&\quad& C_{\alpha\beta}/2
\end{eqnarray*}
\item[\ 4.] If the word $D$ consists of just one letter $\alpha$ then $\varepsilon_{(^\alpha_0)(^\beta_0)}=C_{\alpha\beta}/2$, $\varepsilon_{(^\alpha_1)(^\beta_0)}=-C_{\alpha\beta}/2$.   If the word $D$ consists of just one letter $\bar{\alpha}$ then $\varepsilon_{(^\alpha_0)(^\beta_0)}=-C_{\alpha\beta}/2$, $\varepsilon_{(^\alpha_1)(^\beta_0)}=C_{\alpha\beta}/2$. 
\end{trivlist}

\vskip 3mm
{\bf Example.} For the word $D=\alpha\bar{\beta}\bar{\alpha}\bar{\alpha}\beta$ considered above one can easily compute that all nonvanishing elements of $\varepsilon$ are given by:

$$
\widehat{\varepsilon}_{(^\alpha_0)(^\alpha_1)}=
-\widehat{\varepsilon}_{(^\alpha_1)(^\alpha_2)}=
-\widehat{\varepsilon}_{(^\alpha_2)(^\alpha_3)}=\widehat{C}_{\alpha \alpha}/2
$$
$$
-\widehat{\varepsilon}_{(^\beta_0)(^\beta_1)}=
\widehat{\varepsilon}_{(^\beta_1)(^\beta_2)}=\widehat{C}_{\beta \beta}/2
$$
$$
\widehat{\varepsilon}_{(^\alpha_0)(^\beta_0)}=
\widehat{\varepsilon}_{(^\alpha_3)(^\beta_2)}=\widehat{C}_{\alpha \beta}/2
$$
$$
\widehat{\varepsilon}_{(^\alpha_0)(^\gamma_0)}=
-\widehat{\varepsilon}_{(^\alpha_1)(^\gamma_0)}/2=
\widehat{\varepsilon}_{(^\alpha_3)(^\gamma_0)}=\widehat{C}_{\alpha \gamma}/2
$$
$$
-\widehat{\varepsilon}_{(^\beta_0)(^\gamma_0)}=
\widehat{\varepsilon}_{(^\beta_1)(^\gamma_0)}/2=
-\widehat{\varepsilon}_{(^\beta_2)(^\gamma_0)}=\widehat{C}_{\beta \gamma}/2
$$
\begin{proposition} \label{2.14.05.1} The definitions \ref{2.13.05.2} and 
\ref{9.23.04.1} are equivalent. 
\end{proposition} 

\vskip 3mm
{\bf Proof.} The property 4. of the 
matrix $\varepsilon$ tells us that the two definitions coincide for the elementary 
seeds $\mathbf J(\alpha)$. So it remains to check that the seed $\mathbf J(D)$ 
is the amalgamated product. 

\vskip 3mm
\paragraph{3.4 A map to the group.}
Recall the torus  ${\cal X}_{\mathbf J(D)}= {\mathbb G}_m^{J(D)}$ and the natural coordinates $\{x^\alpha_i\}$ on it. Let us define the map $ev: {\cal X}_{\mathbf J(D)}\rightarrow G$. In order to do it we are going to construct a sequence of group elements, each of which is either a constant or depends on just one coordinate of ${\cal X}_{\mathbf J(D)} $. The product of the elements of the sequence would give the desired map.

Let $f_\alpha,h_\alpha,e_\alpha$ be Chevalley generators of the Lie algebra $\mathfrak g$ of $G$. They are defined up to an action of the Cartan group $H$ of $G$. 
Let $\{h^\alpha\}$ be another basis of the Cartan subalgebra defined by the property:
$$
[h^\alpha,e_\beta]=\delta^\alpha_\beta e_j,\quad [h^\alpha,f_\beta]
=-\delta^\alpha_\beta f_\beta.
$$
This basis is related to the basis 
$\{h_\alpha\}$ via the  Cartan matrix:
$\sum_\beta C_{\alpha \beta} h^\beta=h_\alpha$.

Recall the lattice $X_*(H)$ of homomorphisms (cocharacters)
 ${\Bbb G}_m \to H$. The elements $h_\alpha$ and $h^\alpha$ give rise to 
cocharacters $H_\alpha, H^\alpha \in X_*(H)$, 
called the coroot and the coweight corresponding to the simple root $\alpha$:
$$
H_\alpha: {\Bbb G}_m \to H, \quad dH_\alpha(1) = h_\alpha, \qquad 
H^\alpha: {\Bbb G}_m \to H, \quad dH^\alpha(1) = h^\alpha 
$$
 One has $H^\alpha(x) = \exp(\log(x)h^\alpha)$. 
Let us introduce the group elements ${\mathbf E}^\alpha=\exp e_\alpha$, ${\mathbf F}^\alpha=\exp f_\alpha $.

Replace the letters in $D$ by the group elements using the rule $\alpha \rightarrow {\mathbf E}^\alpha, \bar{\alpha}\rightarrow {\mathbf F}^\alpha$. For $\alpha \in \Pi$ and for any $(^\alpha_i)\in J(D)$ insert $H^\alpha(x^\alpha_i)$ somewhere between the corresponding walls. The choice in placing every $H$ is unessential since it commutes with all $E$'s and $F$'s unless they are marked by the same root.

In other words the sequence of group elements is defined by the following requirements:
\begin{itemize}
\item The sequence of ${\mathbf E}$'s and ${\mathbf F}$'s reproduce the sequence of the letters in the word $D$.
\item Any $H$ depends on its own variable $x^\alpha_i$.
\item There is at least one ${\mathbf E}^\alpha$ or ${\mathbf F}^\alpha$ between any two $H^\alpha$'s.
\item The number of $H^\alpha$'s is equal to the total number of ${\mathbf E}^\alpha$'s and
${\mathbf F}^\alpha$'s plus one.
\end{itemize}

\vskip 3mm
{\bf Example:} The word $\alpha\bar{\beta}\bar{\alpha}\bar{\alpha}\beta$ is mapped by $ev$ to
$$
H^\alpha(x^\alpha_0)H^\beta(x^\beta_0){\mathbf E}^\alpha H^\alpha(x^\alpha_1) {\mathbf F}^\beta
H^\beta(x^\beta_1) {\mathbf F}^\alpha  H^\alpha(x^\alpha_2){\mathbf F}^\alpha {\mathbf E}^\beta
H^\alpha(x^\alpha_3)H^\beta(x^\beta_2)H^\gamma(x^\gamma_0)
$$

\vskip 3mm
\paragraph{3.5. The key properties of the spaces ${\cal X}_B$.} 
We are going to show that the association of the seed ${\cal X}$-torus 
to the words is compatible with natural operations on the words.
\begin{theorem} \label{9.23.04.2} 
Let $A,B$ be arbitrary words and $\alpha,\beta\in \Pi$. Then 
there are the following rational maps commuting with the map $ev$:
\begin{enumerate}
\item ${\cal X}_{\mathbf J(A\bar{\alpha}\beta B)} \rightarrow {\cal
  X}_{\mathbf J(A\beta\bar{\alpha} B)}$ if $\alpha \neq \beta$.
\item ${\cal X}_{\mathbf J(A\alpha\beta B)} \rightarrow {\cal X}_{\mathbf J(A\beta\alpha B)}$ and
${\cal X}_{\mathbf J(A\bar{\alpha}\bar{\beta} B)} \rightarrow {\cal X}_{\mathbf
J(A\bar{\beta}\bar{\alpha} B)}$\\ if $C_{\alpha\beta}=0$.
\item ${\cal X}_{\mathbf J(A\alpha\beta\alpha B)} \rightarrow 
{\cal X}_{\mathbf J(A\beta\alpha\beta
B)}$ and ${\cal X}_{\mathbf J(A\bar{\alpha}\bar{\beta}\bar{\alpha} B)} 
\rightarrow {\cal X}_{\mathbf
J(A\bar{\beta}\bar{\alpha}\bar{\beta} B)}$ \\ if $C_{\alpha\beta}=-1$, 
$C_{\beta\alpha} = -1$.
\item ${\cal X}_{\mathbf J(A\alpha\beta\alpha\beta B)} \rightarrow {\cal
  X}_{\mathbf J(A\beta\alpha\beta\alpha B)}$ and ${\cal X}_{\mathbf
J(A\bar{\alpha}\bar{\beta}\bar{\alpha}\bar{\beta}  B)} \rightarrow {\cal X}_{\mathbf
J(A\bar{\beta}\bar{\alpha}\bar{\beta}\bar{\alpha} B)}$ \\ if $C_{\alpha\beta}=-2,  
C_{\beta\alpha} = -1$.
\item  ${\cal X}_{\mathbf J(A\alpha\beta\alpha\beta\alpha\beta B)} \rightarrow {\cal X}_{\mathbf
J(A\beta\alpha\beta\alpha\beta\alpha B)}$ and ${\cal X}_{\mathbf
J(A\bar{\alpha}\bar{\beta}\bar{\alpha}\bar{\beta}\bar{\alpha}\bar{\beta} B)} \rightarrow {\cal
X}_{\mathbf J(A\bar{\beta}\bar{\alpha}\bar{\beta}\bar{\alpha}\bar{\beta}\bar{\alpha} B)}$\\ if
$C_{\alpha\beta} = -3, C_{\beta\alpha} = -1$.
\item ${\cal X}_{\mathbf J(A\alpha\alpha B)} \rightarrow {\cal X}_{\mathbf J(A\alpha B)}$
and ${\cal X}_{\mathbf J(A\bar{\alpha}\bar{\alpha}B)} \rightarrow {\cal X}_{\mathbf
J(A\bar{\alpha}B)}$.
\item ${\cal X}_{\mathbf J(A)}\times {\cal X}_{\mathbf J(B)}\to {\cal X}_{\mathbf J(AB)}$.
\end{enumerate}
The maps 1,2 are isomorphisms. The maps 3,4,5 are cluster transformations. (They are composions of $1$, $3$ and at least $10$ mutations, respectively). The map 6 is a composition of a cluster transformation and a projection along the coordinate axis. The map 7 is an amalgamated product. The map $ev$ is a Poisson map. The maps 1 through 7 are also Poisson maps.
\end{theorem}

\vskip 3mm
{\bf Remark.} It is not true that a mutation of a cluster seed $\mathbf J(D)$ is always a cluster seed corresponding to another word of the semigroup. For example, let $\Pi=\{\gamma,\Delta,\eta\}$ be the root system of type $A_3$ 
 with $C_{\eta\gamma}=0$. (It is handy to use the capital letter $\Delta$ to distinguish the root which 
plays a specual role below). Then $\mu_{(^\gamma_1)} \mathbf J(\gamma\Delta\eta\gamma\Delta\gamma) = \mathbf J(\Delta\gamma\Delta\eta\Delta\gamma)$, but $\mu_{(^\Delta_1)} \mathbf J(\gamma\Delta\eta\gamma\Delta\gamma)$ is a seed which does not correspond to any word.

\vskip 3mm
{\bf Remark.} There is a famous 8-terms relation between the relations in the symmetric group. Namely in the notations of the previous remark we have:
$$
\gamma\Delta\eta\gamma\Delta\gamma=
\gamma\Delta\eta\gamma\Delta\gamma=
\gamma\eta\Delta\eta\gamma\Delta\stackrel{\sim}{=}
\eta\gamma\Delta\gamma\eta\Delta=
\eta\Delta\gamma\Delta\eta\Delta=
\eta\Delta\gamma\eta\Delta\eta\stackrel{\sim}{=}
\eta\Delta\eta\gamma\Delta\eta=
$$
$$
=\Delta\eta\Delta\gamma\Delta\eta=
\Delta\eta\gamma\Delta\gamma\eta\stackrel{\sim}{=}
\Delta\gamma\eta\Delta\eta\gamma=
\Delta\gamma\Delta\eta\Delta\gamma=
\gamma\Delta\gamma\eta\Delta\gamma=
\gamma\Delta\eta\gamma\Delta\gamma
$$
It is equivalent to the relation between mutations:
$$
\mu_{(^\Delta_1)}\mu_{(^\gamma_2)}\mu_{(^\gamma_1)}\mu_{(^\Delta_1)}\mu_{(^\gamma_2)}\mu_{(^\gamma_1)}\mu_{(^\Delta_1)}\mu_{(^\gamma_2)}=id
$$
It is an easy exercise to show that this relation is a corrolary of the properties 5 and 6 of mutations. Thus the 8-term relation can be reduced to pentagons in the cluster setting.

\vskip 3mm
{\bf Proof}. The proof of the last two statements follow immediately from the 
rest of the theorem.  To prove the rest of the theorem we define the map 7 as the amalgamation map between the corresponding ${\cal X}$-varieties. Then it is sufficient to construct the other maps for the shortest word where the maps are defined, and then extend them using the multiplicativity property 7 to the general case. 
The claim that the evaluation map is Poisson will be proved in the subsection 3.8. 

It is handy to recall an explicit description of the amalgamation map 
${\cal X}_{\mathbf J(A)}\times{\cal X}_{\mathbf J(B)}\rightarrow {\cal X}_{\mathbf J(AB)}$. 
Let $\{x_i^\alpha\}$, $\{y_i^\alpha\}$ and $\{z_i^\alpha\}$ be the coordinates on ${\mathcal X}_{\mathbf J(A)}$, ${\mathcal X}_{\mathbf J(B)}$ and ${\mathcal X}_{\mathbf J(AB)}$, respectively. Then the map is given by the formula:
$$
z^\alpha_i = \left\{\begin{array}{lcl}
             x^\alpha_i &\mbox{if} & i<n^\alpha(A)\\
             x^\alpha_{n^\alpha(A)}y^\alpha_0 &\mbox{if} & i=n^\alpha(A)\\
             y^\alpha_{i+n^\alpha(A)} &\mbox{if} & i>n^\alpha(A)
             \end{array}
             \right.
$$
The crucial point is that, just by the construction, this map is compatible with the evaluation map  $ev$ to the group. 
 
The maps 1-6 and their properties are deduced from the Proposition \ref{9.25.04.1} below. Some formulas of this proposition are equivalent to results available in the literature (\cite{L1}, \cite{BZ}, \cite{FZ}), but stated there in a different form. Our goal is to make apparent their cluster nature, i.e. to show that they transform as the
${\cal X}$-coordinates for a cluster varieties. In the non simply-laced cases these transformations are presented as
compositions of several cluster transformations. The very existence of such presentations is a key new result.

\begin{proposition} \label{9.25.04.1} There are the following  identities between the  generators ${\mathbf E}^\alpha, H^\alpha(x), {\mathbf E}^\alpha$:
\begin{trivlist}

\item[\ $\underline{\alpha\alpha\to\alpha}:$] $ {\mathbf E}^\alpha H^\alpha(x){\mathbf E}^\alpha 
= H^\alpha(1+x) {\mathbf E}^\alpha H^\alpha(1+x^{-1})^{-1} $\vspace{5mm}

\item[\ $\underline{\alpha\beta\to\beta\alpha}:$] If $C_{\alpha\beta}=0$. then 
$\mathbf E^\alpha \mathbf E^\beta = \mathbf E^\beta \mathbf E^\alpha$,\vspace{5mm}

\item[\ $\underline{\alpha\beta\alpha\to\beta\alpha\beta}:$] 
If $C_{\alpha\beta}=-1$, then 
$\mathbf E^\alpha H^\alpha(x) \mathbf E^\beta \mathbf E^\alpha = \\
= H^\alpha(1+x)H^\beta(1+x^{-1})^{-1} \mathbf E^\beta H^\beta(x)^{-1}
\mathbf E^\alpha \mathbf E^\beta H^\alpha(1+x^{-1})^{-1} H^\beta(1+x)$ \vspace{5mm}

\item[\ $\underline{\alpha\beta\alpha\beta\to\beta\alpha\beta\alpha}:$] If
  $C_{\alpha\beta}=-2, C_{\beta\alpha}=-1$,  then 
$\mathbf E^\alpha \mathbf E^\beta H^\alpha(x) H^\beta(y) \mathbf E^\alpha \mathbf E^\beta  =\\=H^\beta(a')H^\alpha(b')\mathbf E^\beta \mathbf E^\alpha H^\beta(y')H^\alpha(x') \mathbf E^\beta \mathbf E^\alpha H^\beta(q')H^\alpha(p')$, 

\noindent where $a',b',x',y',p',q'$ are rational functions of $x$ and $y$ given by (\ref{b2-mutation})
\vspace{5mm}

\item[\ $\underline{\alpha\beta\alpha\beta\alpha\beta\to\beta\alpha\beta\alpha\beta\alpha}:$] If $C_{\alpha\beta}=-3, C_{\beta\alpha}=-1$, then \begin{equation}\label{g2-seqz}
\begin{array}{c}
\mathbf E^\alpha H^\alpha(x)\mathbf E^\beta H^\beta(y)\mathbf E^\alpha H^\alpha(z)\mathbf E^\beta H^\beta(w)\mathbf E^\alpha \mathbf E^\beta=\\
=
H^\beta(a')H^\alpha(b')\mathbf E^\beta H^\beta(y')\mathbf E^\alpha H^\alpha(x')\mathbf E^\beta H^\beta(w')\mathbf E^\alpha H^\alpha(z')\mathbf E^\beta H^\beta(q')\mathbf E^\alpha H^\alpha(p')
\end{array}
\end{equation}
\noindent where $a',b',x',y',z',w',p',q'$ are rational functions of $x,y,z,w$ given by (\ref{g2-mutation})
\end{trivlist}

Further, one has 
\begin{trivlist}
\item[\ $\underline{\bar{\alpha}\alpha\to\alpha\bar{\alpha}}:$] $\mathbf F^\alpha H^\alpha(x) 
\mathbf E^\alpha =\\= \left(\prod\limits_{\beta \neq  \alpha} H^\beta(1+x)^{-C_{\alpha\beta}}\right) 
H^\alpha(1+x^{-1})^{-1}  \mathbf E^\alpha H^\alpha(x^{-1})\mathbf F^\alpha H^\alpha(1+x^{-1})^{-1}$

\item[\ $\underline{\bar{\alpha}\beta\to\beta\bar{\alpha}}:$] $\mathbf F^\alpha \mathbf E^\beta =  \mathbf E^\beta \mathbf F^\alpha$, if $\alpha \neq \beta$.
\end{trivlist}
\end{proposition} 

Applying the antiautomorphism of ${\mathfrak g}$ which acts as  the identity on the Cartan subalgebra, and 
interchanges ${\mathbf F}^\alpha$ and ${\mathbf E}^\alpha$, we obtain similar formulas 
for ${\bar{\alpha}\bar{\alpha}\to\bar{\alpha}}$, 
${\bar{\alpha}\bar{\beta}\to\bar{\beta}\bar{\alpha}}$, ${\bar{\alpha}\bar{\beta}\bar{\alpha}\to\bar{\beta}\bar{\alpha}\bar{\beta}}$ and 
${\bar{\alpha}\bar{\beta}\bar{\alpha}\bar{\beta}\to\bar{\beta}\bar{\alpha}\bar{\beta}\bar{\alpha}}$. For example
 if $C_{\alpha\beta}=C_{\beta\alpha}=-1$, then 
$$
\mathbf F^\alpha H^\alpha(x) \mathbf F^\beta \mathbf F^\alpha = 
= H^\alpha(1+x^{-1})^{-1}H^\beta(1+x) \mathbf F^\beta H^\alpha(x)^{-1}
\mathbf F^\alpha \mathbf F^\beta H^\alpha(1+x) H^\beta(1+x^{-1})^{-1} 
$$

\vskip 3mm
{\bf Proposition \ref{9.25.04.1} implies Theorem \ref{9.23.04.2} minus ``$5$'' and ``the map $ev$ is Poisson'' parts}. 
The map 1 is a corollary of the obvious property $\underline{\bar{\alpha}\beta\to\beta\bar{\alpha}}$. The maps 2 follow from $\underline{\alpha\beta\to\beta\alpha}$ and $\underline{\bar{\alpha}\bar{\beta}\to\bar{\beta}\bar{\alpha}}$. 
The maps $6$ follow from $\underline{\alpha\alpha\to\alpha}$ and $\underline{\bar{\alpha}\bar{\alpha}\to\bar{\alpha}}$. 
Each of these maps is obviously a composition of a mutation $\mu_{(^\alpha_{n^\alpha(A)+1})}$, or its $\bar \alpha$ version, and the 
projection along the corresponding coordinate. The maps $3$ follow from  $\underline{\alpha\beta\alpha\to\beta\alpha\beta}$ 
and $\underline{\bar{\alpha}\bar{\beta}\bar{\alpha}\to\bar{\beta}\bar{\alpha}\bar{\beta}}$ and 
are given by the mutations  
$\mu_{(^\alpha_{n^\alpha(A)+1})}$, or its $\bar \alpha$ version. 
The maps $4$ follow form 
$\underline{\alpha\beta\alpha\beta\to\beta\alpha\beta\alpha}$ 
and $\underline{\bar{\alpha}\bar{\beta}\bar{\alpha}\bar{\beta}
\to\bar{\beta}\bar{\alpha}\bar{\beta}\bar{\alpha}}$ 
and can be easily shown to be given by composition 
of three mutations: $\mu_{(^\alpha_{n^\alpha(A)+1})}$, 
$\mu_{(^\beta_{n^\beta(A)+1})}$ and $\mu_{(^\alpha_{n^\alpha(A)+1})}$, or their bar counterparts. 
A conceptual proof explaining this  see in 
subsection 3.7 below. However it is not at all clear from very complicated 
fromulas (\ref{g2-mutation}) why the 
map $5$ in Theorem \ref{9.23.04.2} is a composition of mutations.

\vskip 3mm
{\bf Proof of Proposition \ref{9.25.04.1}}. It will occupy the end of this and the next two subsections, and will be combined with the proof of the part ``$5$'' of Theorem \ref{9.23.04.2}, as well as a more conceptual proof of the part ``$4$''.
 
Recall that there is a  
transposition antiautomorphism which interchanges $\mathbf E^\alpha$ and $\mathbf F^\alpha$ and does not change $H^\alpha$. Thus the formulae  $\underline{\bar{\alpha}\bar{\alpha}\to\bar{\alpha}}$, $\underline{\bar{\alpha}\bar{\beta}\to\bar{\beta}\bar{\alpha}}$, $\underline{\bar{\alpha}\bar{\beta}\bar{\alpha}\to\bar{\beta}\bar{\alpha}\bar{\beta}}$ and $\underline{\bar{\alpha}\bar{\beta}\bar{\alpha}\bar{\beta}\to\bar{\alpha}\bar{\beta}\bar{\alpha}\bar{\beta}}$ follow from the respective formulae for positive roots.  

 Let us introduce a more traditional generator $H_\alpha(x)=\exp(\log(x)h_\alpha)$.
\vskip 3mm
\noindent{$\underline{\alpha\alpha\to\alpha}$}. It is easy to show using computations with $2\times2$ matrices that
$$
\mathbf E^\alpha H_\alpha(t)\mathbf E^\alpha = H_\alpha(1+t^2)^{1/2} 
\mathbf E^\alpha H_\alpha(1+t^{-2})^{-1/2}
$$
Substituting  
$H_\alpha(t)=\prod_\beta H^\beta(t)^{C_{\alpha \beta}}$ and taking
into account that $\mathbf E^\alpha$ and $\mathbf H^\beta$ commute when $\beta\neq 
\alpha$,  
$C_{\alpha \alpha}=2$,
 and making a
substitution 
$x= t^2$ one gets the identity {$\underline{\alpha\alpha\to\alpha}$}.

\vskip 3mm
\noindent {$\underline{\bar{\alpha}\beta\to\beta\bar {\alpha}}$}. 
The proof is similar to the previous one. It is based on 
the easily verifiable identity for
$2\times2$ matrices
$$
\mathbf F^\alpha H_\alpha(t) \mathbf E^\alpha = H_\alpha(1+t^{-2})^{-1/2}\mathbf E^\alpha H_\alpha(t)^{-1} \mathbf F^\alpha
H_\alpha(1+t^{-2})^{-1/2}
$$

\vskip 3mm
\noindent{$\underline{\alpha\beta\alpha\to\beta\alpha\beta}$.} This identity can be easily derived from the well known identity, which is sufficient to check for $SL_3$: 
$$e^{ae_\alpha}e^{be_\beta}e^{ce_\alpha}=e^{\frac{bc}{a+c}e_\beta}e^{(a+c)e_\alpha}e^{\frac{ab}{a+c}e_\beta}.$$ 
Taking into account that 
\begin{equation} \label{newold}
H^\alpha(a)\mathbf E^\alpha H^\alpha(a)^{-1}=e^{a e_\alpha}
\end{equation} 
for any $\alpha$ and making the substitution $a/c \rightarrow x$ one obtains the desired identity.

\vskip 3mm
\paragraph{3.6  Cluster folding.} 
We start from recalling the notion of the {\em folding} of root systems. Let $\Pi'$ and $\Pi$ be two sets fo simple roots corresponding to the root systems with the  Cartan matrices $C'$ and $C$, respectively. A surjective map $\pi:\Pi'\rightarrow \Pi$ is called folding if it satisfies the following properties:
\vskip 2mm
1. $C_{\alpha'\beta'}=0$ if $\pi(\alpha')=\pi(\beta')$, and $\alpha'\not = \beta'$. 

2. $C_{\alpha,\beta}=\sum_{\alpha' \in \pi^{-1}(\alpha)}C_{\alpha'\beta'}$ if $\pi(\beta')=\beta$. 
\vskip 2mm

A folding induces an embedding (in the inverse direction) of the corresponding Lie algebras denoted by $\pi^*$ and given by:
$$
\pi^*(h_{\alpha})=\sum_{\alpha' \in \pi^{-1}(\alpha)}h'_\alpha, \quad \pi^*(e_{\alpha})=\sum_{\alpha' \in \pi^{-1}(\alpha)}e'_\alpha, \quad \pi^*(e_{-\alpha})=\sum_{\alpha' \in \pi^{-1}(\alpha)}e'_{-\alpha},
$$
where $\{h'_{\alpha'},e'_{\alpha'},e'_{-\alpha'}\}$ are the standard Chevalley generators of the Lie algebra $\mathfrak g'$ corresponding to the Cartan matrix $C'$. 

A folding induces also maps between the corresponding Weyl groups, braid semigroups and braid groups, and Hecke semigroups, given by
$$
\pi^*(\alpha)=\prod_{\alpha' \in \pi^{-1}(\alpha)}\alpha', \quad \pi^*(\bar{\alpha})=\prod_{\alpha' \in \pi^{-1}(\alpha)}\bar{\alpha}'.
$$ 
(The order of the product does not matter since according to the property 1 the factors commute.)

The main feature of the folding is that it gives embeddings of non simply-laced Lie algebras and groups to the simply-laced ones. Namely, $B_n$ is a folding of $D_n$, $C_n$ is a folding of $A_{2n-1}$, $F_4$ is a folding of $E_6$ and $G_2$ is a folding of both $B_3$ and $D_4$. In these cases, the folding is provided by the action of a subgroup $\Gamma$ of the automorphism group  of the Dynkin diagram: one has $\Pi':= \Pi/\Gamma$, and the folding is the quotient map 
$\Pi \to \Pi/\Gamma$.  On the level of Dynkin diagrams a folding corresponds just to the folding of the corresponding graph, thus explaining the origin of the name. 

\vskip 3mm
Following \cite{FG2}, we define a {\it folding of a cluster seed}:

\begin{definition} \label{2.14.04.1g} A {\em folding} $\pi$ of a cluster seed $\mathbf J'=(J',J'_0,\varepsilon',d')$ to a cluster seed $\mathbf J=(J,J_0,\varepsilon,d)$ is a surjective map $\pi:J'\rightarrow J$ satisfying the following conditions:
\vskip 2mm
0. $\pi(J'_0) = J'_0$, $\pi(J'-J'_0) = J-J_0$.

1. $\varepsilon'_{i'j'}=0$ if $\pi(i')=\pi(j')$.

2. $\varepsilon_{ij}= \sum_{i' \in \pi^{-1}i}\varepsilon'_{i'j'}$, and all summands in this sum have the same signs or vansih. 

\end{definition}
\vskip 3mm 
A folding induces a map of the corresponding cluster tori $\pi^*: \mathcal X_{\mathbf J} \rightarrow \mathcal X_{\mathbf J'}$ by the formula $(\pi^*)^*x_{i'}=x_{\pi(i')}$. The main feature of this map is that it commutes with mutations in the following sense. To formulate it we need the following simple but basic fact.
\begin{lemma} \label{2.14.04.1}
If $\varepsilon_{k k'}=0$, then mutations in the directions $k$ and $k'$.
commute. 
\end{lemma}
\vskip 3mm
Let $\mu_k:\mathbf J \rightarrow \mathbf I$ be a mutation, $\pi:\mathbf J'\rightarrow \mathbf J$ be a folding, and  $\mathbf I':=(\prod_{k'\in \pi^{-1}(k)}\mu_{k'})\mathbf J'$. We define a map  $\pi_k: \mathbf I' \rightarrow \mathbf I$ as a composition 
$$
\pi_k=\mu_k\pi(\prod_{k'\in \pi^{-1}(k)}\mu_{k'})^{-1}
$$ 
The last factor in this formula is well defined since the mutations $\mu_{k'}$ commute thanks to Lemma \ref{2.14.04.1} and the condition 1 of Definition \ref{2.14.04.1g}. The map $\pi_k$ is not always a folding (the condition 2 may not be satisfied), but if it is then,  of course, $\mu_k\pi = \pi_k\prod_{k' \in \pi^{-1}k}\mu_{k'}$; furthermore, on the level of $\mathcal X$-tori we have $\pi^*\mu_k=\prod_{k' \in \pi^{-1}k}\mu_{k'}\pi^*$.

We would like to notice that the folding map is not a Poisson map. However it sends symplectic leaves 
to the symplectic leaves, and, being restricted to a symplectic leave, multiplies the symplectic structure there 
by a constant. 
\vskip 3mm
The two foldings: of the Cartan matrices and of the cluster seeds, are closely related. Namely, let 
$\pi:\Pi'\rightarrow \Pi$ be a folding of the Cartan matrices. Denote by $\mathfrak W$ 
(respectively $\mathfrak W'$) the free seimgroup generated by $\Pi$ and $-\Pi$ 
(respectively by $\Pi'$ and $-\Pi'$). 
 Let $D \in \mathfrak W$, and let $D'=\pi^*(D)$ be the image of $D$ in the semigroup $\mathfrak W'$. 
The proof of the following proposition is rather straitforward, and thus left to the reader. 

\begin{proposition} In the above notation, there is a  natural map $\pi:\mathbf J(D')
\rightarrow \mathbf J(D)$, which is a folding of cluster seeds. Moreover the map of the corresponding seed $\mathcal X$-tori commutes with the evaluation map $ev$ to the respective Lie groups.   
\end{proposition}

\vskip 3mm
\paragraph{3.7 A proof of the parts ``$4$'' and ``$5$'' of Theorem \ref{9.23.04.2}.}  Let us prove first the formula $\underline{\alpha\beta\alpha\beta\to\beta\alpha\beta\alpha}$. We need the folding map of the simple roots $\Pi'$ of the Lie group $A_3$ to the simple roots $\Pi$ of the Lie group $B_2$. Let $\Pi'=\{\gamma,\Delta,\eta\}$, $\Pi=\{\alpha,\beta\}$, $\pi(\gamma)=\pi(\eta)=\alpha$, $\pi(\Delta)=\beta$,
$$
\left(
\begin{array}{ccc}
      C'_{\gamma\gamma}&C'_{\gamma\Delta}&C'_{\gamma\eta}\\
      C'_{\Delta\gamma}&C'_{\Delta\Delta}&C'_{\Delta\eta}\\
      C'_{\eta\gamma}&C'_{\eta\Delta}&C'_{\eta\eta}
\end{array}
\right)
=
\left(
\begin{array}{rrr}
      2&-1&0\\
      -1&2&-1\\
      0&-1&2
\end{array}
\right),
\qquad
\left(
\begin{array}{cc}
       C_{\alpha\alpha}&C_{\alpha\beta}\\
       C_{\beta\alpha}&C_{\beta\beta}
\end{array}
\right)
=
\left(
\begin{array}{rr}
       2&-2\\
       -1&2
\end{array}
\right)
$$
There is the following sequence of relations in the braid group of type $A_3$:
\begin{equation}\label{b2seq}
\pi^*(\alpha\beta\alpha\beta)=
\gamma\eta\Delta\gamma\eta\Delta\stackrel{\sim}{=}
\eta\gamma\Delta\gamma\eta\Delta=\eta\Delta\gamma\Delta\eta\Delta = 
\eta\Delta\gamma\eta\Delta\eta\stackrel{\sim}{=}\eta\Delta\eta\gamma\Delta\eta=
$$ 
$$=\Delta\eta\Delta\gamma\Delta\eta=
\Delta\eta\gamma\Delta\gamma\eta\stackrel{\sim}{=}
\Delta\gamma\eta\Delta\gamma\eta=\pi^*(\beta\alpha\beta\alpha),
\end{equation}
which just shows that $\pi^*$ is a semigroup homomorphism. Here $\stackrel{\sim}{=}$ stands for the elementary transformations provided by the relation $\gamma\eta = \eta\gamma$.

Hence we have

\begin{equation}\label{b2xseq}
\pi^*(\mathbf E^\alpha \mathbf E^\beta H^\alpha(x) H^\beta(y) \mathbf E^\alpha \mathbf E^\beta)=
\end{equation}
$$
\mathbf E^\gamma \mathbf E^\eta \mathbf E^\Delta H^\gamma(x)H^\eta(x) H^\Delta(y) \mathbf E^\gamma \mathbf E^\eta \mathbf E^\Delta=\mathbf E^\eta \mathbf E^\gamma \mathbf E^\Delta H^\gamma(x)H^\eta(x) H^\Delta(y) \mathbf E^\gamma \mathbf E^\eta \mathbf E^\Delta=
$$
$$=H^\gamma(1+x)H^\Delta(1+x^{-1})^{-1}\mathbf E^\eta \mathbf E^\Delta H^\Delta(x^{-1}) \mathbf E^\gamma \mathbf E^\Delta H^\Delta(z+xz) H^\eta(y) \mathbf E^\eta \mathbf E^\Delta H^\gamma(1+x^{-1})^{-1}=\cdots
$$
$$
\cdots= H^\eta(a')H^\gamma(a')H^\Delta(b')\mathbf E^\Delta \mathbf E^\gamma \mathbf E^\eta H^\gamma(x')H^\eta(x') H^\Delta(y') \mathbf E^\Delta \mathbf E^\gamma \mathbf E^\eta H^\eta(p')H^\gamma(p')H^\Delta(q') = 
$$
$$
=\pi^*(H^\beta(b')H^\alpha(a')\mathbf E^\beta \mathbf E^\alpha H^\beta(y')H^\alpha(x') \mathbf E^\beta \mathbf E^\alpha H^\beta(q')H^\alpha(p')),
$$
where the ellipsis $\cdots$ means repeated application of the formula $\underline{\alpha\beta\alpha\to\beta\alpha\beta}$ corresponding to the sequence of mutations  $\mu_1^\gamma\mu_1^\eta\mu_1^\Delta\mu_1^\gamma $, and 
\begin{equation}\label{b2-mutation}
\begin{array}{ll}
a'=\frac{1+x+2xy+xy^2}{1+x+xy},& \quad b'=\frac{xy^2}{1+x+2xy+xy^2},\\
p'= 1+x+xy,& \quad q'= \frac{x(1+x+2xy+x^2y)}{(1+x+xy)^2} ,\\
x'=\frac{y}{1+x+2xy+xy^2},& \quad y'=\frac{(1+x+xy)^2}{xy^2}
\end{array}
\end{equation}
This proves the $\underline{\alpha\beta\alpha\beta\to\beta\alpha\beta\alpha}$ claim of Proposition \ref{9.25.04.1}. 
Further, from this we easyly get, by adding pairs of elements of the Cartan group on both sides of (\ref{b2xseq}), a birational 
transformation
$$
\Psi_{B_2}: \Q(a, b, p, q, x, y) \lra \Q(a'', b'', p'', q'', x'', y'') 
$$
which reduces to (\ref{b2-mutation}) when $a=b=p=q=1$, and is determined by the formula
\begin{equation}\label{b2xseq1}
\pi^*(H^\alpha(a)H^\beta(b)\mathbf E^\alpha \mathbf E^\beta H^\alpha(x) H^\beta(y) \mathbf E^\alpha \mathbf E^\beta H^\alpha(p)H^\beta(q))= 
\end{equation}
$$
\pi^*(H^\beta(b'')H^\alpha(a'')\mathbf E^\beta \mathbf E^\alpha H^\beta(y'')H^\alpha(x'') \mathbf E^\beta \mathbf E^\alpha H^\beta(q'')H^\alpha(p'')).
$$
To prove that the map 4 from Theorem 3.5 is a cluster transformation of the original $\mathcal X$-torus, we need to show that there exists a sequence of mutations of the cluster seed $\mathbf J(\alpha\beta\alpha\beta)$ whose product is equal to the transformation $\Psi_{B_2}$.

Consider  a cluster transformation $L_B: 
{\cal X}_{\alpha \beta\alpha \beta} \lra 
{\cal X}_{\beta\alpha \beta\alpha}$ given as a composition of three mutations: 
$$
L_B:= \mu_{(^\alpha_1)} \mu_{(^\beta_1)} \mu_{(^\alpha_1)}
$$
Let us show that the map $\Psi_{B_2}$ is equal to 
the cluster transformation $L_B$.

The $4$ nontrivial transformations in (\ref{b2seq}) 
give rise to a cluster transformation 
$$
L_A: {\cal X}_{\gamma\eta\Delta \gamma\eta\Delta } \lra 
{\cal X}_{\Delta \gamma\eta\Delta \gamma\eta}, \quad L_A= \mu_{(^\gamma_1)}\mu_{(^\eta_1)}\mu_{(^\Delta_1)}
\mu_{(^\gamma_1)}
$$
given by composition of the corresponding sequence of $4$ mutations. 
There is a  diagram 
\begin{equation} \label{2.16.05.5}
\begin{array}{ccc}
{\cal X}_{\alpha \beta\alpha \beta}&\hra&
{\cal X}_{\gamma\eta\Delta \gamma\eta\Delta }\\
L_B\downarrow &&\downarrow L_A\\
{\cal X}_{\beta\alpha \beta \alpha }&\hra&{\cal X}_{\Delta \gamma\eta\Delta \gamma\eta}
\end{array}
\end{equation}
where the horizontal arrows are the folding embeddings. 
\begin{lemma}\label{2.16.05.6e}
The diagram (\ref{2.16.05.5}) is  commutative. 
\end{lemma}

\begin{figure}[ht]
\centerline{\epsfbox{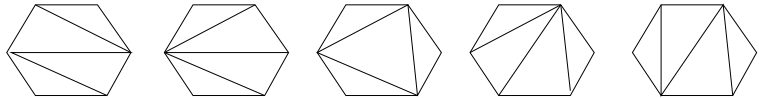 }}
\caption{The sequence of $4$ mutations corresponding to $L_A$}
\label{cox1}
\end{figure}

\begin{figure}[ht]
\centerline{\epsfbox{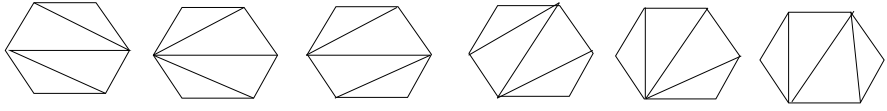}}
\caption{The sequence of $5$ mutations corresponding to $ \widehat L_B $}
\label{cox2}
\end{figure}

{\bf Proof}. Consider the following cluster transformation, which is 
the image under the folding embedding of $L_B$: 
\begin{equation} \label{2.16.05.6er}
\widehat L_B: {\cal X}_{\gamma\eta\Delta \gamma\eta\Delta } 
\lra {\cal X}_{\delta \gamma\eta\delta \gamma\eta}, \qquad 
\widehat L_B:= \mu_{(^{\eta}_1)} \mu_{(^{\gamma}_1)} 
\mu_{(^\Delta_1)} 
\mu_{(^{\eta}_1)} \mu_{(^{\gamma}_1)}
\end{equation}
It evidently makes the following diagram commutative 
\begin{equation} \label{2.16.05.6}
\begin{array}{ccc}
{\cal X}_{\alpha \beta\alpha \beta}&\hra&
{\cal X}_{\gamma\eta\Delta \gamma\eta\Delta }\\
 L_B\downarrow &&\downarrow \widehat L_B\\
{\cal X}_{\beta\alpha \beta \alpha }&\hra&{\cal X}_{\Delta \gamma\eta\Delta  \gamma
\eta}
\end{array}
\end{equation}
It remains to show that the cluster transformations $\widehat L_B$ and 
$L_A$ are equal.
It can be done by computing the effect of the action of the latter sequence of mutations on  the $\mathcal X$-torus, and checking that it coincides with the transformation $\Psi_{B_2}$. Another way is to use explicitly the pentagon relations. For the conaisseurs of the cluster varieties we give another proof just by drawning pictures. The seed $\mathbf J(\gamma\eta\Delta\gamma\eta\Delta)$ is of the finite type $A_3$, it has only a finite number of different seeds parametrized by triangulations of a hexagon. (See Appendix 2, where we discuss this model and the isomorphism of the corresponding cluster variety with the configuration space of $6$ points in ${\Bbb P}^1$ and of $3$ flags in $PGL_4$). In this framework mutations correspond to removing an edge of the triangulation and replacing it by another diagonal of the arising quadrilateral. Thus the first sequence of mutations corresponds to the sequence of triangualtions shown on Fig.\ref{cox1}, 

\noindent while the second one is shown on Fig. \ref{cox2}.
The Lemma, and hence the part 4 of Theorem \ref{9.23.04.2} are proved.

Now let us proceed to the proof of the formula $\underline{\alpha\beta\alpha\beta\alpha\beta\to\beta\alpha\beta\alpha\beta\alpha}$. Consider the folding map $\pi$ of the set of simple roots $\Pi'$ of the Lie algebra $D_4$ to the set of simple roots $\Pi$ of the Lie algebra $G_2$. Let $\Pi'=\{\gamma,\eta,\rho,\Delta\}$, $\Pi=\{\alpha,\beta\}$, $\pi(\gamma)=\pi(\eta)=\pi(\rho)=\alpha$, $\pi(\Delta)=\beta$.

$$
\left(
\begin{array}{cccc}
      C'_{\gamma\gamma}&C'_{\gamma\eta}&C'_{\gamma\rho}&C'_{\gamma\Delta}\\
      C'_{\eta  \gamma}&C'_{\eta  \eta}&C'_{\eta  \rho}&C'_{\eta  \Delta}\\
      C'_{\rho  \gamma}&C'_{\rho  \eta}&C'_{\rho  \rho}&C'_{\rho  \Delta}\\
      C'_{\Delta\gamma}&C'_{\Delta\eta}&C'_{\Delta\rho}&C'_{\Delta\Delta}
\end{array}
\right)
=
\left(
\begin{array}{rrrr}
       2& 0& 0&-1\\
       0& 2& 0&-1\\
       0& 0& 2&-1\\
      -1&-1&-1& 2
\end{array}
\right),
\qquad
\left(
\begin{array}{cc}
       C_{\alpha\alpha}&C_{\alpha\beta}\\
       C_{\beta\alpha}&C_{\beta\beta}
\end{array}
\right)
=
\left(
\begin{array}{rr}
       2&-3\\
       -1&2
\end{array}
\right)
$$

There is the following sequence of relations in the braid group of type $D_4$:
\begin{equation}\label{g2-xseq}
\begin{array}{c}
\pi^*(\alpha\beta\alpha\beta\alpha\beta)=\\
\gamma\eta  \rho  \Delta\gamma\eta  \rho  \Delta\gamma\eta  \rho  \Delta\stackrel{\sim}{=}
\gamma\eta  \rho  \Delta\rho  \eta  \gamma\Delta\gamma\eta  \rho  \Delta=
\gamma\eta  \Delta\rho  \Delta\eta  \gamma\Delta\gamma\eta  \rho  \Delta=
\gamma\eta  \Delta\rho  \Delta\eta  \Delta\gamma\Delta\eta  \rho  \Delta=\\
=
\gamma\eta  \Delta\rho  \eta  \Delta\eta  \gamma\Delta\eta  \rho  \Delta\stackrel{\sim}{=}
\gamma\eta  \Delta\eta  \rho  \Delta\gamma\eta  \Delta\eta  \rho  \Delta=
\gamma\Delta\eta  \Delta\rho  \Delta\gamma\eta  \Delta\eta  \rho  \Delta=
\gamma\Delta\eta  \Delta\rho  \Delta\gamma\Delta\eta  \Delta\rho  \Delta=\\
=
\gamma\Delta\eta  \Delta\rho  \Delta\gamma\Delta\eta  \rho  \Delta\rho  =
\gamma\Delta\eta  \rho  \Delta\rho  \gamma\Delta\eta  \rho  \Delta\rho  \stackrel{\sim}{=}
\gamma\Delta\eta  \rho  \Delta\gamma\rho  \Delta\rho  \eta  \Delta\rho  =
\gamma\Delta\eta  \rho  \Delta\gamma\Delta\rho  \Delta\eta  \Delta\rho  =\\
=
\gamma\Delta\eta  \rho  \gamma\Delta\gamma\rho  \Delta\eta  \Delta\rho  \stackrel{\sim}{=}
\gamma\Delta\gamma\eta  \rho  \Delta\rho  \gamma\Delta\eta  \Delta\rho  =
\Delta\gamma\Delta\eta  \rho  \Delta\rho  \gamma\Delta\eta  \Delta\rho  =
\Delta\gamma\Delta\eta  \Delta\rho  \Delta\gamma\Delta\eta  \Delta\rho  =\\
=
\Delta\gamma\eta  \Delta\eta  \rho  \Delta\gamma\Delta\eta  \Delta\rho  =
\Delta\gamma\eta  \Delta\eta  \rho  \Delta\gamma\eta  \Delta\eta  \rho  \stackrel{\sim}{=}
\Delta\gamma\eta  \Delta\rho  \eta  \Delta\eta  \gamma\Delta\eta  \rho  =
\Delta\gamma\eta  \Delta\rho  \Delta\eta  \Delta\gamma\Delta\eta  \rho  =\\
=
\Delta\gamma\eta  \Delta\rho  \Delta\eta  \gamma\Delta\gamma\eta  \rho  =
\Delta\gamma\eta  \rho  \Delta\rho  \eta  \gamma\Delta\gamma\eta  \rho  \stackrel{\sim}{=}
\Delta\gamma\eta  \rho  \Delta\gamma \eta \rho  \Delta\gamma\eta  \rho  =
\pi^*(\beta\alpha\beta\alpha\beta\alpha)
\end{array}
\end{equation}
It just shows that $\pi^*$ is a homomorphism of semigroups. 
Here $\stackrel{\sim}{=}$ stands for the equalities which follow from the commutativity relations $\gamma\eta=\eta\gamma$, $\gamma\rho=\rho\gamma$ and $\eta\rho=\rho\eta$.

A computation similar to (\ref{b2xseq}) is too long to write it here; it was done using a computer. The result is:
\begin{equation}\label{g2-seq}
\begin{array}{c}
\pi^*(\mathbf E^\alpha H^\alpha(x)\mathbf E^\beta H^\beta(y)\mathbf E^\alpha H^\alpha(z)\mathbf E^\beta H^\beta(w)\mathbf E^\alpha \mathbf E^\beta)=\\
=
\pi^*(H^\beta(a')H^\alpha(b')\mathbf E^\beta H^\beta(y')\mathbf E^\alpha H^\alpha(x')\mathbf E^\beta H^\beta(w')\mathbf E^\alpha H^\alpha(z')\mathbf E^\beta H^\beta(q')\mathbf E^\alpha H^\alpha(p'))
\end{array}
\end{equation}
where
\begin{equation}\label{g2-mutation}
\begin{array}{ll}
a'=\frac{xR_2}{R_3},& \quad b'=R_3,\\
p'= \frac{xyz^2w}{R_1},& \quad q'= \frac{R_1^3}{R_4} ,\\
x'=\frac{zR_1R_3}{R_4},& \quad y'=\frac{yR_4}{R_2^3}\\
z'=\frac{R_4}{xyz^2wR_2},& \quad w'= \frac{wR_2^3}{R_3^3}
\end{array}
\end{equation}
and
\begin{eqnarray*}
R_1&=&xyz^2w+1+x+yx+2yxz+yxz^2\\
R_2&=&y^2x^2z^3w+y^2x^2z^3+3y^2x^2z^2+3yx^2z+2yxz+3y^2x^2z+1+2x+2yx+\\& &+x^2+2yx^2+y^2x^2\\
R_3&=&3x+3x^2+3yx^2+3yx^2z+1+y^2x^3z^3w+y^2x^3z^3+3y^2x^3z^2+3yx^3z+3y^2x^3z+\\& &+x^3+2yx^3+y^2x^3\\
R_4&=&1+3x+3y^2x^3z^4+12y^2x^3z+6yx^3z+18y^2x^3z^2+12y^2x^3z^3+x^3+2y^2x^3z^3w+\\& &+3y^2x^2z^4w+3y^2x^2z^3w+3yx+3y^2x^3z^4w+6yxz+3x^2+6yx^2+12yx^2z+3yxz^2+\\& &+3y^2x^3+3yx^3+3yx^3z^2+2y^3x^3z^6w+6y^3x^3z^5w+6y^3x^3z^4w+y^3x^3+20y^3x^3z^3+\\& &+6y^3x^3z^5+6y^3x^3z+15y^3x^3z^4+y^3x^3z^6+15y^3x^3z^2+2y^3x^3z^3w+y^3x^3z^6w^2+\\& &+6yx^2z^2+12y^2x^2z+18y^2x^2z^2+12y^2x^2z^3+3y^2x^2z^4+3y^2x^2.
\end{eqnarray*}
This proves the $\underline{\alpha\beta\alpha\beta\alpha\beta\to\beta\alpha\beta\alpha\beta\alpha}$ claim of Proposition \ref{9.25.04.1}. 

Just like in the $B_2$ case, to prove that the map 5 from the theorem 3.5 is a cluster transformation, we need to show that there exists a cluster transformation of the seed $\mathbf J(\alpha\beta\alpha\beta\alpha\beta)$ which, being transformed by the folding map to the $D_4$ set-up,  equals to the cluster transformation 
encoded in the sequence (\ref{g2-xseq}), and given explicitly as the left hand side in the formula 
(\ref{opopo}) below. To do this it is sufficient to show the following equality between two sequences of mutations:
\begin{equation} \label{opopo}
\mu_{(^\rho_1)}\mu_{(^\rho_2)}\mu_{(^\Delta_1)}\mu_{(^\eta_2)}\mu_{(^\eta_1)}\mu_{(^\rho_1)}\mu_{(^\gamma_1)}
\mu_{(^\gamma_2)}\mu_{(^\rho_2)}\mu_{(^\rho_1)}\mu_{(^\Delta_2)}\mu_{(^\eta_2)}\mu_{(^\eta_1)}\mu_{(^\Delta_1)}\mu_{(^\gamma_2)}\mu_{(^\rho_1)}=
\end{equation}
$$
\mu_{(^\Delta_2)}\mu_{(^\rho_1)}\mu_{(^\eta_1)}\mu_{(^\gamma_1)}\mu_{(^\Delta_1)}\mu_{(^\Delta_2)}\mu_{(^\rho_2)}\mu_{(^\eta_2)}\mu_{(^\gamma_2)}\mu_{(^\Delta_2)}\mu_{(^\rho_1)}\mu_{(^\eta_1)}\mu_{(^\gamma_1)}\mu_{(^\rho_2)}
\mu_{(^\eta_2)}\mu_{(^\gamma_2)}\mu_{(^\Delta_1)}\mu_{(^\Delta_2)}
$$
Here the second composition is the image under the  folding map of the sequence of mutations 
$\mu_{(^\beta_2)}\mu_{(^\alpha_1)}\mu_{(^\beta_1)}\mu_{(^\beta_2)}\mu_{(^\alpha_2)}\mu_{(^\beta_2)}\mu_{(^\alpha_1)}
\mu_{(^\alpha_2)}\mu_{(^\beta_1)}\mu_{(^\beta_2)}$.

This was done by an explicit calculation of the action of the latter sequences of mutations on the ${\cal X}$-coordinates, performed by a computer, which showed that it coincides with the transformation given by (\ref{g2-seq})-(\ref{g2-mutation}). It would be interesting to find a proof which relates one of 
 the  sequences of mutations in the Lie group of type $D_4$ to the other by using the pentagon relations. 

\vskip 3mm
\paragraph{3.8 The evaluation map $ev$ 
is a Poisson map.} To prove this claim, it is sufficient to prove it in the simplest case: 

\begin{proposition} \label{alphapoisson} 
Let $\alpha$ be a simple root. Then the evaluation map ${ev}: 
\mathcal{X}_{\mathbf J(\alpha)} \hookrightarrow G$ is a Poisson immersion. So 
its image is a Poisson subvariety of $G$, and the induced Poisson structure on $\mathcal{X}_{\mathbf J(\alpha)}$ coincides with the one (\ref{Poisson_bracket}) 
for the matrix $\varepsilon$ given by (\ref{Elementary_epsilon}). 
\end{proposition}
\vskip 3mm

We deduce the general claim from the proposition by induction, using the following four facts: 
the multiplication map $G \times G \to G$ is a Poisson map 
for the standard Poisson structure on $G$, the evaluation map commutes 
with the multiplication, i.e. the 
diagram (\ref{mevcomm1}) is commutative,  and the left vertical map in that diagram 
 is a Poisson map, and a dominant map,  i.e. its image is dence, and thus the induced map 
of functions is injective. 

\vskip 3mm
{\bf Proof of Proposition \ref{alphapoisson}}. The evaluation map is obviously immersion in our case. 
Let us recall the standard Poisson structure on $G$. Let $R \subset \mathfrak h^{\ast}$ be the set of roots of the Lie algebra $\mathfrak{g}$ with a Cartan subalgebra $\mathfrak{h}$, and let $R_+ \subset R$ be the subset of positive roots. The root decomposition of $\mathfrak{g}$ reads as $\mathfrak{g} = \oplus_{\beta \in R}\  \mathfrak{g}_{\beta} \oplus \mathfrak{} \mathfrak{h}$. Let $\{e_{\alpha} \in \mathfrak{g}_{\alpha}| \alpha \in R \}$ be a set of root vectors normalized so that $[e_{\alpha}, [ e_{- \alpha}, e_{\alpha} ] ] = 2 e_{\alpha}$. The vectors $e_{\pm \alpha}$, $\alpha \in R_+$, are defined by this condition uniquely up to rescaling $e_{\pm \alpha} \to 
\lambda^{\pm 1}_{\alpha}e_{\pm \alpha}$. Let $r \in \mathfrak{g} \otimes \mathfrak{g}$ be the standard $r$-matrix:
\begin{equation}
 r = \sum_{\alpha \in R_+} d^\alpha e_{\alpha} \wedge e_{- \alpha} = \sum_{\alpha \in R_+} d^\alpha (e_{\alpha} \otimes e_{- \alpha} - e_{- \alpha} \otimes e_{\alpha} ) ,
\label{r-matrix}
\end{equation}
where as above $d^\alpha=\frac{(\alpha,\alpha)}{2}$.
Observe that $r$ does not depend on the choice of the vectors $e_{\alpha}$ satisfying the above normalization. 
The Poisson bracket on $G$ is given by a bivector field $P = r^L - r^R $, where $r^L$ (resp. $r^R$) is the right-invariant (resp. left-invariant) 
bivector field on $G$ which equals $r$ at the identity of $G$. If we identify the tangent space to $G$ at a point $g \in G$ with $\mathfrak{g}$ using the right shift, the value $P(g)$ of the bivector field  $P$ at $g$ is  
$P(g) = r - \mbox{Ad}_g r$.  
We apply this formula in the special case when  
\begin{equation} \label{g}
g = (\prod_{\beta \in \Pi } H^{\beta} ( x_0^{\beta} )) \mathbf{E}^\alpha H^\alpha (x_1^\alpha )
\end{equation} 
To make the computation we shall use the following formulae:
\begin{equation}
 \mbox{Ad}_{H^{\alpha}(x)} r = r, \quad
 \mbox{Ad}_{{H^\alpha}(x)} e_\alpha \wedge h_\alpha = x e_\alpha \wedge h_\alpha, \quad 
\mbox{Ad}_{\mathbf E^{\alpha}} r = r + d^\alpha e_\alpha \wedge h_\alpha \label{lemma1},
\end{equation}
The first two are obvious; a proof of the third one see below. Using them, one easily derives 
\begin{equation} \label{biv}
P(g) = r - \mbox{Ad}_g r = d^\alpha x^\alpha_0 h_\alpha \wedge e_\alpha.
\end{equation}
 So to find the Poisson bracket induced on $\mathcal{X}_\alpha$ we need to compute the right-invariant vector fields on $\mathcal{X}_{\alpha}$ corresponding to $h_{\alpha}$ and $e_\alpha$. Obviously, $h_\alpha$ gives rise to a vector field $\sum_{\beta} C_{\alpha \beta } x^\beta_0 \frac{\partial}{\partial x^\beta_0}$. 
The following computaiton, where $g$ is from (\ref{g}),  shows that 
$e_{\alpha}$ gives rise to 
$\frac{\partial}{\partial x_0^\alpha} - \frac{x_1^\alpha}{x_0^\alpha} 
\frac{\partial}{\partial x_1^\alpha}$:
$$ e_\alpha g = \frac{d}{dt} \exp t e_\alpha g |_{t = 0} = \frac{d}{dt}H^\alpha (t)\mathbf E^\alpha H^\alpha (t^{-1})g|_{t = 0} =
$$
$$
 = \frac{d}{dt}(\prod_{\beta \in \Pi - \{ \alpha \}} H^\beta(x_0^\beta)) 
H^{\alpha}(t)\mathbf E^\alpha H^\alpha(t^{- 1}) H^\alpha(x_0^\alpha )\mathbf E^\alpha H^\alpha (x_1^\alpha)|_{t = 0} =
$$
$$
= (\prod_{\beta \in \Pi - \{ \alpha \}} H^\beta (x_0^\beta))\frac{d}{dt} H^\alpha (x_0^\alpha + t)\mathbf E^\alpha H^\alpha (x_1^\alpha - t x_1^\alpha / x_0^\alpha)|_{t = 0} = 
\left(\frac{\partial}{\partial x_0^\alpha} - \frac{x_1^\alpha}{x_0^\alpha}\frac{\partial}{\partial x_1^{\alpha}}\right) g.
$$
Here we have used the formula $\underline{\alpha \alpha \rightarrow \alpha}$ from Proposition \ref{9.25.04.1}.

Subsitiuting these expressions for the vector fields into  (\ref{biv}),
we get 
$$
P =-
\sum_\beta d^\alpha C_{\alpha \beta} x^\beta_0 x_0^\alpha \left(\frac{\partial}{\partial x_0^\alpha} - \frac{x_1^{\alpha}}{x_0^{\alpha}} \frac{\partial}{\partial x_1^\alpha} \right) \wedge \frac{\partial}{\partial x^\beta_0}= 
$$
\begin{equation} \label{5.12.05.1} 
=-\widehat{C_{}}_{\alpha \beta } x_0^\alpha x^\beta_0  \frac{\partial}{\partial x_0^\alpha} \wedge \frac{\partial}{\partial x_0^\beta} + \widehat{C}_{\alpha \beta}x_1^\alpha x^\beta_0 \frac{\partial}{\partial x_1^\alpha} \wedge \frac{\partial}{\partial x_0^\beta}, \end{equation} 
which  coincides with the expression given by (\ref{Poisson_bracket}) and (\ref{Elementary_epsilon}).

{\it Proof of the third formula in (\ref{lemma1})}.  
Let $R (\alpha) \subset R$ be the
 root system for the Dynkin diagram obtained from the initial one 
by deleting  
the vertex corresponding to the simple positive root $\alpha$. Let  $R_+(\alpha)$ be 
the set of its positive roots. Let $p$ be  the projection of $\mathfrak{h}^{\ast}$ onto its quotient by the subspace spanned by  $\alpha$. 
Then $p(R) = R(\alpha)$. 
So we can rewrite the root decomposition as:
\begin{equation}
\mathfrak{g} = \oplus_{\beta \in R(\alpha)}\left(\oplus_{\gamma \in p^{-1} (\beta)}
\mathfrak{g}_{\gamma}\right) 
\oplus \mathfrak{h} (\alpha) \oplus \mathfrak i_{\alpha}(\mathfrak{sl}_2), \label{alpha-decomp}
\end{equation}
where  $i_{\alpha}(\mathfrak{sl}_2)$ is the $\mathfrak{sl}(2)$-subalgebra spanned by $e_\alpha, e_{-\alpha}$, $h_\alpha = [ e_\alpha, e_{-\alpha}]$ and $\mathfrak h (\alpha) \subset \mathfrak h$ is the kernel of $\alpha$. 
Observe that the summands are $i_{\alpha}(\mathfrak{sl}_2)$-invariant.

Let us consider the quadratic Casimir:
$$
t = \sum_{\beta \in R}d^\beta e_{\beta} \otimes e_{-\beta} + \frac12 \sum_{\beta \in \Pi} d^{\beta} h_{\beta} \otimes h^{\beta}\in \mathfrak g \otimes \mathfrak g.$$
It can be rewritten as
\begin{equation}\label{casimir-decomp}
t = t_0 + \sum_{\beta \in R(\alpha)} t_{\beta},
\end{equation}
where $$
t_{\beta} = \sum_{\gamma \in p^{- 1}(\beta )} d^\gamma
e_{\gamma} \otimes e_{-\gamma}, \qquad 
t_0 = d^\alpha (e_{\alpha} \otimes e_{-\alpha}+e_{-\alpha} \otimes e_{\alpha}) +\frac12 \sum_{\beta \in \Pi} d^{\beta} h_{\beta}\otimes h^{\beta}.
$$
Every term of the expression (\ref{casimir-decomp})  is $i_{\alpha}(\mathfrak{sl}_2)$-invariant. 
The  $r$-matrix (\ref{r-matrix}) is decomposed in the same way:
$$r = r_0 + \sum_{\beta \in R_+(\alpha)} t_{\beta} 
- \sum_{\beta \in R_- ( \alpha )} t_{\beta},$$
where $t_{\beta}$ are as above and $r_0 = d^{\alpha}e_\alpha \wedge e_{-\alpha}.$ All the terms but the first one $r_0$ are $i_{\alpha}(\mathfrak{sl}_2)$-invariant and thus $\mbox{Ad}_{\mathbf E^\alpha}(r) = r-r_0+\mbox{Ad}_{\mathbf{E}^{\alpha}}(r_0) = r-r_0 + \exp (\mbox{ad}_{e_\alpha}) (e_{\alpha} \wedge e_{-\alpha}) = r + d^\alpha e_\alpha \wedge h_\alpha$.
We proved the third formula in (\ref{lemma1}), and hence the Proposition \ref{alphapoisson}. As it was explained above, this implies that $ev$ is a Poisson map in general. 

Therefore, we completed the proof of Theorem \ref{9.23.04.2}. 

\vskip 3mm
\paragraph{3.9 Duality conjectures and canonical bases.} 
Let $w_0$ be the longest element of the Weyl group $W_G$ of $G$. 
Let ${\cal X}_{G, w_0}$ be the corresponding cluster ${\cal X}$-variety. 
Recall that, given a seed ${\bf I}$, we defined in Section 2 of \cite{FG2} a positive space 
${\cal A}_{|{\bf I}|}$  assigned to it. 
Let $G^L$ be the Langlands dual group for $G$. 
Applying this construction to the seed corresponding to the element $w_0$ 
in $G^L$, we arrive at 
a positive space ${\mathcal A}_{G^L, w_0}$. Recall that for any semifield $\mathbb F$ and a positive space ${\cal X}$ 
there is a set ${\mathcal X}(\mathbb F)$ of $\mathbb F$-points of ${\mathcal X}$, (loc. cit.). 
Recall the tropical semifield $\mathbb Z^t$: it is the set $\mathbb Z$ with the following semifield operations: the  multiplication and division are given by the usual addition and subtraction in $\Z$, and the semifield addition is given by taking the maximum.  Let ${\cal A}_{G^L, w_0}(\Z^t)$ be the set of $\Z^t$-points of the positive space ${\cal A}_{G^L, w_0}$.  

Then, according to the duality conjecture from Section 5 of \cite{FG2}, there should exist a basis in the algebra $\Z[{\cal X}_{G, w_0}]$ of regular functions on the variety ${\mathcal X}_{G, w_0}$, parametrised by the set ${\cal A}_{G^L, w_0}(\Z^t)$. Let us explain how it should be related to the (dual) canonical basis of Lusztig \cite{L2}. 

The dual canonical basis is a basis of regular functions on the Borel subgroup $B$ of $G$. Since ${\cal X}_{G, w_0}$ is birationally equivalent to $B$, the regular functions on the former are rational, but not necessarily regular, functions on $B$. Let us say that an element of  our conjectural basis in $\mathbb Z[{\cal X}_{G, w_0}]$ is {\it regular} if it provides a regular function on $B$. 

\begin{conjecture} \label{5-3.05.1}
The regular elements of the conjectural basis in $\mathbb Z[{\mathcal X}_{G, w_0}]$ form a basis of the space of regular functions on $B$. Moreover, it coincides with Lusztig's dual canonical basis  on $B$. 
\end{conjecture} 
\vskip 3mm

Let us try to determine when a rational function $F$ on  ${\cal X}_{G, w_0}$ is regular on $B$. Observe that $B = HU$, where $U$ is the maximal unipotent in $B$. Apparently $F$ is regular on $H$. So it remains to determine when it is regular on $U$. 

Pick a reduced decomposition of $w_0$. Let $(t_{i_1}, ..., t_{i_N})$, where $N = {\rm dim}U = l(w_0)$,  be the corresponding Lusztig's coordinates on $B$ (\cite{L1}), and $(x_{i_1}, ..., x_{i_N+r})$ the corresponding cluster ${\cal X}$-coordinates on ${\cal X}_{G, w_0}$. Choose coordinates $(h_1, .., h_r)$ on $H$. Evidently  the $t$-coordinates are related to the $(x,h)$-coordinates by monomial transformations, and vise versa. In particular any  $F \in \mathbb Z[{\cal X}_{G, w_0}]$ is a Laurent polynomial in $(t_{i_1}, ..., t_{i_N})$. 

\begin{lemma} \label{5-4.05.1}
An $F \in \Z[{\cal X}_{G, w_0}]$ is regular on $B$ if and only if for any reduced decomposition of $w_0$ it is a polynomial in the corresponding coordinates 
$(t_{i_1}, ..., t_{i_N})$. 
\end{lemma}

\vskip 3mm
{\bf Proof}.  The ``only if'' part is clear. Let us check the opposite. 
The subvarieties $ev({\cal X}_{G, w})$ in $G$, where $w\in {W}_G$, are of codimension $l(w_0) - l(w)$, and it follows from the Bruhat decomposition that the 
complement to their union is of codimension $\geq 2$. Any irreducible component of divisors $ev({\cal X}_{G, w})$ is given by the equation $t_{i_k}=0$ for certain reduced decomposition of $w_0$ and certain $k$. The lemma follows.

\begin{figure}[ht]
\centerline{\epsfbox{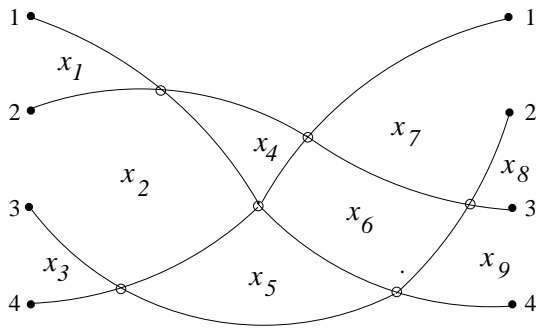}}
\caption{The wiring diagram and ${\cal X}$-coordinates for the word 
$\sigma_3\sigma_1 \sigma_2\sigma_1\sigma_3\sigma_2$. }
\label{cox6}
\end{figure}

\begin{figure}[ht]
\centerline{\epsfbox{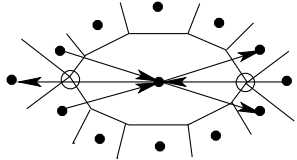}}
\caption{Reading the Poisson tensor from a  wiring diagram.}
\label{cox8}
\end{figure}

\vskip 3mm
\paragraph{3.10 Examples for $PGL_m$.} Below we show how to visualize, in the case of $PGL_m$,
 the combinatorics of the cluster 
${\cal X}$-coordinates 
by a wiring diagram. The wiring diagram languige  is well known \cite{BFZ96}. 
Our goal is to show how it works for the ${\cal X}$-coordinates. 
\vskip 3mm
{\bf 1.} In the case of $PGL_m$ one can 
visualize the combinatorics of the cluster 
${\cal X}$-coordinates related to a word 
by a wiring diagram. The ${\cal X}$-coordinates are assigned to the 
connected components of the complement to the wiring diagram, except the bottom and top components. 
The frozen ${\cal X}$-coordinates are assigned to the very left and right domains, i.e. to the 
domains which are not completely bounded by wires.  
The word itself is encoded by the wiring diagram as follows: we scan the diagram from the left to the right, and 
assign a generator for each vertex of the diagram: the generator  $\sigma_i$ is assigned to a vertex having  
$i-1$ wires above it. See Figure \ref{cox6}, 
illustrating the situation for the word 
$\sigma_3\sigma_1 \sigma_2\sigma_1\sigma_3\sigma_2$ for $PGL_4$. 
The corresponding parametrization of the Borel subgroup of upper triangular $4\times 4$ matrices 
is given by the product
\begin{equation} \label{wd}
{\bf H}^3(x_3){\bf H}^1(x_1){\bf E}^3{\bf E}^1{\bf H}^2(x_2){\bf H}^1(x_4){\bf H}^3(x_5){\bf E}^2{\bf E}^1{\bf H}^2(x_6){\bf E}^3{\bf H}^1(x_7){\bf E}^2{\bf H}^2(x_8){\bf H}^3(x_9)
\end{equation}
Here ${\bf E}^i$ is the elementary unipotent 
matrix corresponding to the $i$-th simple positive 
root: it has $1$'s on the diagonal, and 
the only non-zero non-diagonal element is
$1$ at the entry $(i, i+1)$. Further, ${\bf H}^j(t) = {\rm diag}(\underbrace{t. \ldots, t}_j, 1, \ldots, 1)$ 
is the diagonal matrix corresponding to the $j$-th simple coroot. The frozen variables are 
$x_1, x_2,. x_3, x_7, x_8. x_9$. To record the expression (\ref{wd}), 
we scan the wiring diagram from the left to the right. The intersection points of wires 
provide the elementary matrices ${\bf E}^i$, while the domains contributes 
the Cartan elements ${\bf H}^j(x)$. Observe that the order of factors in (\ref{wd}) 
is by no means uniquely determined: the Cartan elements commute, and some of them 
commute with some ${\bf E}$'s. The wiring diagram, considered modulo isotopy, 
encodes the element (\ref{wd}) in a more adequate way. 
\vskip 3mm
{\bf  2}. The Poisson structure tensor is encoded by the wiring diagram as follows. 
Take a (connected) domain of the wiring diagram corresponding to a non-frozen coordinate $x_0$. 
It can be rather complicated, sharing boundary with many other domains, see 
Figure \ref{cox8}. However it has two distinguished vertices, 
the very left and right ones, shown by circles. There are at most six 
outside domains sharing one of these two vertices. Let $x_j$ be an ${\cal X}$-coordinate 
for the given wiring diagram. 
The Poisson bracket $\{x_0, x_j\} = \varepsilon_{0j}x_0x_j$ 
is non-zero if and only if $x_j$ is assigned to one of those domains. One has 
$\varepsilon_{0j}= \pm 1$, and the sign is shown by arrows on the picture: 
$\varepsilon_{0j}= 1$ if and only if the arrow goes from $x_0$ to $x_j$.
The Posson bracket between two frozen variables is obtained similarly, 
but the coefficient is divided by $2$. 
 
\begin{figure}[ht]
\centerline{\epsfbox{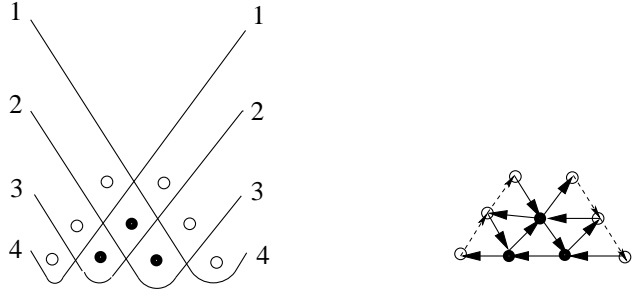}}
\caption{The ${\cal X}$-coordinates for a standard reduced decomposition 
of $w_0$ for $PGL_4$. }
\label{cox7}
\end{figure}

\vskip 3mm
{\bf  3}. Yet another example, corresponding to a ``standard'' reduced decomposition of $w_0$ for $PGL_m$, 
 is given on the left hand side of Figure \ref{cox7}. The non-frozen coordinates, shown by black points, 
give ries to coordinates on $H\backslash B/H$. Observe that there is a canonical birational isomorphism 
between $H\backslash B/H$ and 
the configuration space of triples of flags ${\rm Conf}_3({\cal B})$. 
Using it, one can show that 
the non-frozen  
coordinates in this case are identified with the canonical coordinates on the 
configuration space of triples of flags  for $PGL_m$ introduced in \cite{FG1}, Section 9.  
The right hand side of Figure \ref{cox7} shows the Poisson tensor. 
Its coefficients between the frozen coordinates are shown by dotted arrows. 
\vskip 3mm
\paragraph{3.11. Cluster ${\cal X}$-variety structure of partial flag varieties.}  
Let $P$ be a parabolic subgroup of $G$, and $P=M_PU_P$ its Levi decomposition. 
So $G/P$ is a partial flag variety. Let $w_0^G$ (resp. $w_0^M$)  
be the longest element of the Weyl group of $G$ (resp. $M_P$). Write $w_0^G = w_0^Mw_0^U$.  
Take a reduced decomposition $\widetilde w^U_0$ of $w_0^U$. It gives rise to a coordinate system 
on $G/P$ as follows. Take the seed corresponding to $\widetilde w_0^U$ 
and the corresponding seed ${\cal X}$-torus ${\cal X}_{\widetilde w^U_0}$. 
The frozen part of the torus ${\cal X}_{\widetilde w^U_0}$ 
is a product $H_L \times H_R$ 
of two Cartan subgroups, called the left ($H_L$) and right ($H_R$) frozen Cartan subgroups. 
The canonical projection ${\cal X}_{\widetilde w^U_0} \to G/P$ provides
 a regular open embedding $H_L\backslash 
{\cal X}_{\widetilde w^U_0} \hra U^{\rm opp}_P \hra  G/P$, where $U^{\rm opp}_P$ is the 
subgroup opposite to the unipotent radical $U_P$. It is the coordinate system on $G/P$ 
corresponding to $\widetilde w^U_0$. It follows from Theorem \ref{9.23.04.2} that the collection 
of coordinate systems on $G/P$ corresponding to different reduced decompositions 
of $w_0^U$ provides a set of cluster ${\cal X}$-coordinate systems. 

\begin{figure}[ht]
\centerline{\epsfbox{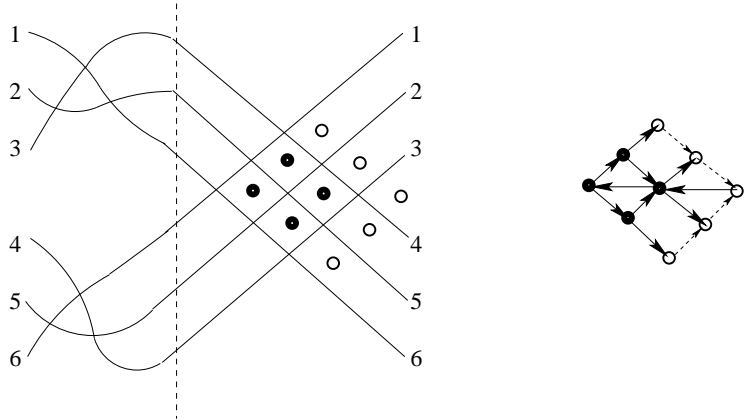}}
\caption{The cluster ${\cal X}$-variety structure of ${\rm Gr}_3(6)$.} 
\label{cox9}
\end{figure}

Observe that $P$ is a Poisson subgroup of $G$, so $G/P$ has a natural Poisson structure. 
It follows from Theorem \ref{9.23.04.2} that it coincides   with the one given by 
the cluster ${\cal X}$-variety structure on $G/P$. 
\vskip 3mm
{\bf Example}. For the Grassmannian ${\rm Gr}_k(n)$ of $k$-planes in 
an $n$-dimensional vector space the above construction provides 
a {\it canonical} coordinate system. Indeed there is 
only one reduced decomposition of $w_0^U$ for the Grassmannian. See Figure \ref{cox9} 
where the case of ${\rm Gr}_3(6)$ is illustrated. The wiring diagram for the Grassmannian is on the right of the dotted vertical line. The frozen variables 
are shown by circles. 
The oriented graph providing the 
Poisson structure tensor is on the right. It is calculated using the recepee from 
Section 3.10. It coincides with the one 
defined in \cite{GSV1}. 

\section{Appendix 1: The braid group of type $G_2$ acts 
on triples of flags of type $G_2$}

In this Section we study the cluster ${\cal X}$-variety corresponding to the 
moduli space ${\rm Conf}_3({\cal B}_{G_2})$ 
of configurations of triples of flags of type $G_2$. The combinatorial structure of 
a cluster ${\cal X}$-variety is reflected in the topology of 
the {\it modular orbifold}, defined in Section 2 of 
\cite{FG2}. Below we recall its definition. Then we 
determine the modular orbifold for the 
moduli space ${\rm Conf}_3({\cal B}_{G_2})$, and compute its fundamental group, 
which turns out to be 
the braid group of type $G_2$. 
 
\vskip 3mm
\paragraph{4.1 The modular orbifold of a cluster ${\cal X}$-variety.} 
It is constructed in three steps: 

1.  We assign to a seed $\mathbf I = (I, I_0, \varepsilon, d)$ 
a simplex $S_{\mathbf I}$ equipped with a bijection of the set of its  
codimension one faces with $I$. It induces a 
bijection between the set of its vertices and $I$: a vertex is labeled by 
the same element as the opposite codimension one face. Recall that an element 
 $k \in I-I_0$ gives rise to a  seed 
mutation  
$\mathbf I \to \mu_k(\mathbf I)$. We  
glue the simplices $S_{\mathbf I}$ and $S_{\mu_k(\mathbf I)}$ along 
their codimension one faces labeled by $k$, matching 
the vertices labeled by the same elements.  
We continue this process by making all possible mutations 
and gluing the corresponding simplices. This 
way we get a simplicial complex $S_{|\mathbf I|}$. 

2. We identify  
simplices corresponding to isomorphic seeds, getting a simplicial complex 
 $\overline S_{|\mathbf I|}$.

3. We remove from  $\overline S_{|\mathbf I|}$ certain faces 
of codimension $\geq 2$ defined as follows. Recall that a seed $\mathbf I$ 
provides a torus ${\cal X}_{\mathbf I}$ with a coordinate system 
$\{x_i^{\mathbf I}\}$. Let $\varphi: \mathbf I \to \mathbf I'$ be a  cluster transformation of seeds. 
We say that it is an  
${\cal X}$-equivalence, if the  induced 
cluster transformation $\varphi_{\cal X}: {\cal X}_{\mathbf I} \to 
{\cal X}_{\mathbf I'}$ is an isomorphism  
respecting the coordinates: $\varphi^*_{\cal X}
x^{\mathbf I'}_{\varphi(i)} =x^{\mathbf I}_{i}$. 
Let  $F$ be a face 
of  $S_{|\mathbf I|}$.  Consider the set of ${\cal X}$-equivalence classes of seeds $\mathbf I'$ such that the simplex $S_{\mathbf I'}$ in 
$S_{|\mathbf I|}$ contains $F$. We say that $F$ is of {\it infinite type} 
if this set is infinite. Removing from 
$\overline S_{|\mathbf I|}$ all faces of infinite type, we get the  
{\it modular orbifold} $M_{|\mathbf I|}$. 

We prove (loc. cit.) that it is indeed 
an orbifold. Its dimension 
is the dimension of the cluster ${\cal X}$-variety minus one. 

\vskip 3mm
\paragraph{4.2 Main results.} In the $\varepsilon$-finite case the {modular orbifold} is glued from 
a finite number of simplices. 
It is non-compact, unless the cluster ${\cal X}$-variety is of finite cluster type. 
In general it can not be compactified by an orbifold, 
but sometimes it can. 
Here is the main result: 

\begin{theorem} \label{3.9.05.1}
a) The cluster ${\cal X}$-variety corresponding to the 
moduli space ${\rm Conf}_3({\cal B}_{G_2})$ is of $\varepsilon$-finite type. 
The number of  non-isomorphic seeds assigned to it is seven. 

b) The corresponding modular complex is homeomorphic to 
$S^3 - L$, where $L$ is a link with two connected components,  and 
$\pi_1(S^3-L)$ is isomorphic to the braid group of type $G_2$.  
\end{theorem}

The {\it mapping class group} of a cluster ${\cal X}$-variety  
was defined in Section 2 of \cite{FG2}. 
It acts by automorphisms of the cluster ${\cal X}$-variety. 
It is always infinite if the cluster structure is of 
$\varepsilon$-finite, but not of finite type. Theorem \ref{3.9.05.1} immediately implies 
the following: 

\begin{corollary}
The mapping class group of this cluster ${\cal X}$-variety corresponding to 
${\rm Conf}_3({\cal B}_{G_2})$ 
is an infinite quotient of 
the braid group of type $G_2$.
\end{corollary}

{\bf Remark}. According to Hypothesis 2.19 from \cite{FG2}, the modular complex of a cluster 
${\cal X}$-variety is the classifying space (in general orbi-space) for the 
corresponding mapping class group.  This plus Theorem \ref{3.9.05.1} imply that the mapping  
class group should be isomorphic to the braid group of type $G_2$. 

\begin{figure}[ht]
\centerline{\epsfbox{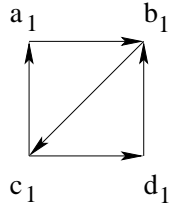}}
\caption{The original seed}
\label{cox13}
\end{figure}

\begin{figure}[ht]
\centerline{\epsfbox{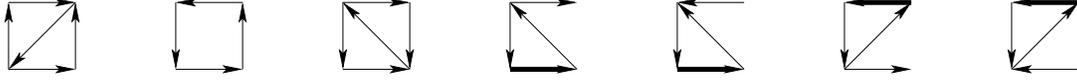}}
\caption{The seven seeds}
\label{cox11}
\end{figure}

\vskip 3mm
\paragraph{4.3 Proof of Theorem \ref{3.9.05.1}.} The cluster structure of the moduli space ${\rm Conf}_3({\cal B}_{G_2})$ can be described by a seed with $4$ vertices. We describe the cluster function by a graph 
 shown on Figure \ref{cox13}. Here two vertices $i,j$ are related by an arrow directed from $i$ to $j$ if and only if 
$\varepsilon_{ij} > 0$. Moreover $\varepsilon_{ij} =1$ in this case. 
The multipliers are equal to $3$ for the top two vertices, and $1$ for the bottom two. 
We claim that mutating this seed we get exactly $7$ different seeds shown on 
Figure \ref{cox11}, where we keep the same convention about the multipliers, while the convention about the arrows is the following: a thin (resp. thick) arrow 
from $i$ to $j$ means that $\varepsilon_{ij}=1$ (resp. $\varepsilon_{ij}=2$). 
To prove this claim we list  below the $14$ pairs of seeds 
from Figure \ref{cox11} related by mutations. 

\begin{figure}[ht]
\centerline{\epsfbox{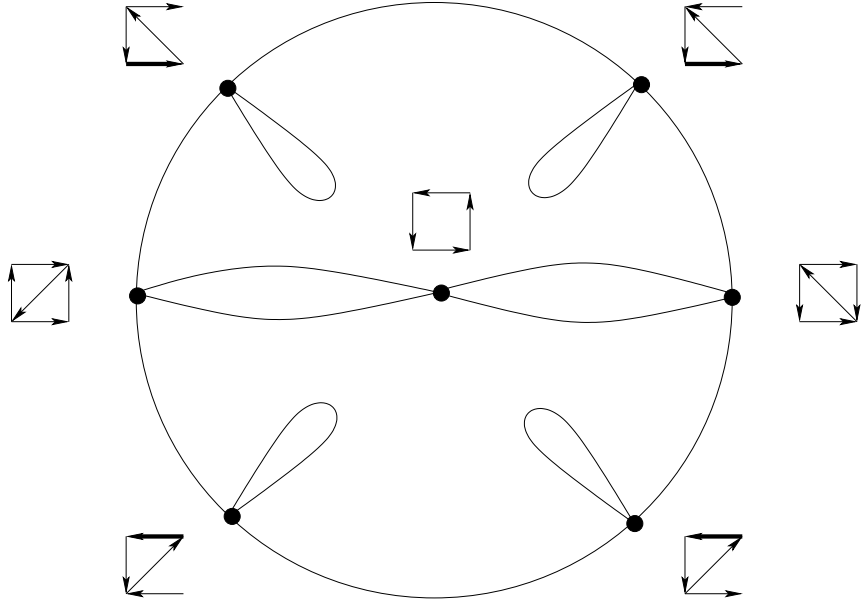}}
\caption{The $1$-skeleton of the dual to the  modular complex}
\label{cox10}
\end{figure}

Let us give an example explaining our notation. We denote by $(a_k, b_k, c_k, d_k)$ the four 
vertices of the seed number $k$ on Figure \ref{cox11} ($k$ counts the seeds from the left to the right). 
The vertices are arranged like on Figure \ref{cox13}: 
$(a_k, b_k)$ are the two top, and $(c_k, d_k)$ the bottom vertices. 
The mutation $\lambda_6$ mutates the seed $(a_2, b_2, c_2, d_2)$ at the vertex $d_2$, producing the 
seed $(a_3, b_3, c_3, d_3)$. Making 
cluster/modular complex, we glue the corresponding two tetrahedra so that the face 
$(a_2, b_2, c_2)$ of the first is glued to the face $(b_3, a_3, d_3)$ of the second, 
matching the $i$-th vertices  of these two triangles. 
We record this information as follows:
$$
\lambda_6 : \begin{pmatrix}a_2 & b_2 &c_2\\b_3 & a_3 &d_3\end{pmatrix} 
$$
Below we list mutations $\lambda_1 , ..., \lambda_{14}$, 
which together with their inverses give us all the mutations: 
$$
 \lambda_1 : \begin{pmatrix}b_1 & c_1 &d_1\\b_2 & c_2 &d_2\end{pmatrix} \quad 
 \lambda_2 : \begin{pmatrix}a_1 & c_1 &d_1\\b_4 & d_4 &c_4\end{pmatrix} \quad 
 \lambda_3 : \begin{pmatrix}a_1 & b_1 &d_1\\a_7 & b_7 &d_7\end{pmatrix} \quad 
 \lambda_4 : \begin{pmatrix}a_1 & b_1 &c_1\\b_2 & a_2 &d_2\end{pmatrix} \quad 
$$
$$
 \lambda_5 : \begin{pmatrix}a_2 & c_2 &d_2\\a_3 & c_3 &d_3\end{pmatrix} \quad 
 \lambda_6 : \begin{pmatrix}a_2 & b_2 &c_2\\b_3 & a_3 &d_3\end{pmatrix} \quad 
 \lambda_7 : \begin{pmatrix}b_3 & c_3 &d_3\\b_5 & d_5 &c_5\end{pmatrix} \quad 
$$
$$
 \lambda_8 :  \begin{pmatrix}a_3 & b_3 &c_3\\a_6 & b_6 &d_6\end{pmatrix} \quad 
 \lambda_9 :  \begin{pmatrix}a_4 & c_4 &d_4\\a_5 & c_5 &d_5\end{pmatrix} \quad 
 \lambda_{10} : \begin{pmatrix}a_4 & b_4 &d_4\\a_4 & b_4 &c_4\end{pmatrix} \quad 
 \lambda_{11} : \begin{pmatrix}a_5 & b_5 &d_5\\a_5 & b_5 &c_5\end{pmatrix} \quad 
$$
$$
 \lambda_{12} : \begin{pmatrix}b_6 & c_6 &d_6\\a_6 & c_6 &d_6\end{pmatrix} \quad 
 \lambda_{13} : \begin{pmatrix}a_6 & b_6 &c_6\\a_7 & b_7 &c_7\end{pmatrix} \quad 
\lambda_{14}: \begin{pmatrix}b_7 & c_7 &d_7\\a_7 & c_7 &d_7\end{pmatrix} \quad 
$$
We present on Figure \ref{cox10} the  $1$-skeleton 
of the simplicial complex dual to the modular complex. It has $7$ vertices corresponding to the 
seven different seeds, and every vertex is connected with the four vertices related 
to it by mutations. 

{\it The gluing data for the edges}. 
Here is how we glue the edges of seven tetrahedra 
making the modular complex. 
$$
\begin{pmatrix}b_1 \\d_1\end{pmatrix} 
\stackrel{\lambda_1}{\lra}\begin{pmatrix}b_2 \\d_2\end{pmatrix}
\stackrel{\lambda_4^{-1}}{\lra}\begin{pmatrix}a_1 \\c_1\end{pmatrix}
\stackrel{\lambda_2}{\lra}\begin{pmatrix}b_4 \\d_4\end{pmatrix}
\stackrel{\lambda_{10}}{\lra}\begin{pmatrix}b_4 \\c_4\end{pmatrix}
\stackrel{\lambda_2^{-1}}{\lra}\begin{pmatrix}a_1 \\d_1\end{pmatrix}
\stackrel{\lambda_3}{\lra}\begin{pmatrix}a_7 \\d_7\end{pmatrix}
\stackrel{\lambda_{14}^{-1}}{\lra}\begin{pmatrix}b_7 \\d_7\end{pmatrix}
\stackrel{\lambda_{3}^{-1}}{\lra}\begin{pmatrix}b_1 \\d_1\end{pmatrix}
$$
$$
\begin{pmatrix}a_2 \\c_2\end{pmatrix} 
\stackrel{\lambda_5}{\lra}\begin{pmatrix}a_3 \\c_3\end{pmatrix}
\stackrel{\lambda_8}{\lra}\begin{pmatrix}a_6 \\d_6\end{pmatrix}
\stackrel{\lambda_{12}^{-1}}{\lra}\begin{pmatrix}b_6 \\d_6\end{pmatrix}
\stackrel{\lambda_{8}^{-1}}{\lra}\begin{pmatrix}b_3 \\c_3\end{pmatrix}
\stackrel{\lambda_7}{\lra}\begin{pmatrix}b_5 \\d_5\end{pmatrix}
\stackrel{\lambda_{11}}{\lra}\begin{pmatrix}b_5 \\c_5\end{pmatrix}
\stackrel{\lambda_{7}^{-1}}{\lra}\begin{pmatrix}b_3 \\d_3\end{pmatrix}
\stackrel{\lambda_{6}^{-1}}{\lra}\begin{pmatrix}a_2 \\c_2\end{pmatrix}
$$
$$
\begin{pmatrix}c_1 \\d_1\end{pmatrix} 
\stackrel{\lambda_1}{\lra}\begin{pmatrix}c_2 \\d_2\end{pmatrix}
\stackrel{\lambda_5}{\lra}\begin{pmatrix}c_3 \\d_3\end{pmatrix}
\stackrel{\lambda_7}{\lra}\begin{pmatrix}d_5 \\c_5\end{pmatrix}
\stackrel{\lambda_9^{-1}}{\lra}\begin{pmatrix}d_4 \\c_4\end{pmatrix}
\stackrel{\lambda_2^{-1}}{\lra}\begin{pmatrix}c_1 \\d_1\end{pmatrix}
$$
$$
\begin{pmatrix}a_1 \\b_1\end{pmatrix} 
\stackrel{\lambda_4}{\lra}\begin{pmatrix}b_2 \\a_2\end{pmatrix}
\stackrel{\lambda_6}{\lra}\begin{pmatrix}a_3 \\b_3\end{pmatrix}
\stackrel{\lambda_8}{\lra}\begin{pmatrix}a_6 \\b_6\end{pmatrix}
\stackrel{\lambda_{13}}{\lra}\begin{pmatrix}a_7 \\b_7\end{pmatrix}
\stackrel{\lambda_3^{-1}}{\lra}\begin{pmatrix}a_1 \\b_1\end{pmatrix}
$$
$$
\begin{pmatrix}b_1 \\c_1\end{pmatrix} 
\stackrel{\lambda_1}{\lra}\begin{pmatrix}b_2 \\c_2\end{pmatrix}
\stackrel{\lambda_6}{\lra}\begin{pmatrix}a_3 \\d_3\end{pmatrix}
\stackrel{\lambda_5^{-1}}{\lra}\begin{pmatrix}a_2 \\d_2\end{pmatrix}
\stackrel{\lambda_4^{-1}}{\lra}\begin{pmatrix}b_1 \\c_1\end{pmatrix}
$$
$$
\begin{pmatrix}a_4 \\c_4\end{pmatrix} 
\stackrel{\lambda_9}{\lra}\begin{pmatrix}a_5 \\c_5\end{pmatrix}
\stackrel{\lambda_{11}^{-1}}{\lra}\begin{pmatrix}a_5 \\d_5\end{pmatrix}
\stackrel{\lambda_9^{-1}}{\lra}\begin{pmatrix}a_4 \\d_4\end{pmatrix}
\stackrel{\lambda_{10}}{\lra}\begin{pmatrix}a_4 \\c_4\end{pmatrix}
$$
$$
\begin{pmatrix}a_6 \\c_6\end{pmatrix} 
\stackrel{\lambda_{12}^{-1}}{\lra}\begin{pmatrix}b_6 \\c_6\end{pmatrix}
\stackrel{\lambda_{13}}{\lra}\begin{pmatrix}b_7 \\c_7\end{pmatrix}
\stackrel{\lambda_{14}}{\lra}\begin{pmatrix}a_7 \\c_7\end{pmatrix}
\stackrel{\lambda_{13}^{-1}}{\lra}\begin{pmatrix}a_6 \\c_6\end{pmatrix}
$$
$$
\begin{pmatrix}c_6 \\d_6\end{pmatrix} 
\stackrel{\lambda_{12}}{\lra}\begin{pmatrix}c_6 \\d_6\end{pmatrix}; \qquad 
\begin{pmatrix}c_7 \\d_7\end{pmatrix}
\stackrel{\lambda_{14}}{\lra}\begin{pmatrix}c_7 \\d_7\end{pmatrix}; \qquad 
\begin{pmatrix}a_4 \\b_4\end{pmatrix}\stackrel{\lambda_{10}}{\lra}
\begin{pmatrix}a_4 \\b_4\end{pmatrix}; \qquad 
\begin{pmatrix}a_5 \\b_5\end{pmatrix} \stackrel{\lambda_{11}}{\lra}\begin{pmatrix}a_5 \\b_5\end{pmatrix}
$$

After the gluing, we get four vertices: 
$$
a_1 = a_2 = a_3 = a_6 = a_7 = b_1 = b_2 = b_3 = b_4 = b_5 = b_6 = b_7,
$$ 
$$
c_1 = c_2 = c_3 = c_4 = c_5 = d_1 = d_2 = d_3 = d_4 = d_5 = d_6 = d_7,
$$ 
$$
a_4 = a_5, \qquad c_6 = c_7
$$  
Thus we have $7$ tetrahedra, $14$ faces, $11$ triangles and $4$ vertices. So the Euler 
characteristic is $0$. 

\begin{figure}[ht]
\centerline{\epsfbox{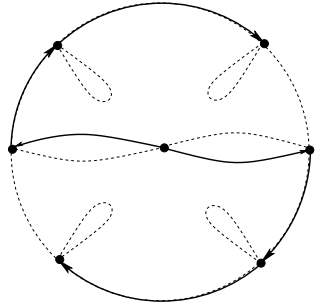}}
\caption{A spanning tree}
\label{cox12}
\end{figure}

Each of the edges $\lambda_i$, $i=1, ..., 14$, of the dual modular complex gives us a generator 
of the fundamental groupoid. Each of the $11$ listed above edges of the original complex 
gives a relation in the fundamental groupoid. 

To compute the fundamental group of the modular complex we use the following algorithm. 
Recall that a spanning tree of a graph is a maximal contractible subgraph of the graph. 
Evidently a  spanning tree contains all vertices of the graph. Let us shrink a spanning tree. Then  
every edge of the graph which does not belong to the spanning tree gives 
rise to a non-trivial loop 
on the quotient. These loops generate the fundamental group of the quotient 
based at the contracted spanning tree.

1. Choose a spanning tree of the dual to the modular complex. 

 2. Then the  fundamental group has the following presentation: 

Generators correspond to the edges of the dual 
modular complex which do not belong to the spanning tree. 
Relations correspond to the edges of the modular complex as follows: Take all triangles 
containing the given edge. 
A co orientation of this edge gives rise to a cyclic order of this set. 
Then the product of the corresponding generators, in an order 
compatible with the cyclic order,   is a relation. 

Let us implement this algorithm. Choose a spanning tree shown by bold arcs on Figure \ref{cox12}.
Then there are $8$ generators, corresponding to mutations 
$\lambda_{3}, \lambda_{4}, \lambda_{5}, \lambda_{7}, 
\lambda_{10}, \lambda_{11}, \lambda_{12}, \lambda_{14}$. The relations can be read off the gluing 
data of the edges:
$$
\lambda_{4}^{-1}\lambda_{10}\lambda_{3}\lambda_{14}^{-1}\lambda_{3}^{-1} = 1, \qquad 
\lambda_{5}\lambda_{12}^{-1}\lambda_{7}\lambda_{11}\lambda_{7}^{-1} = 1, \qquad 
\lambda_{5}\lambda_{7} = 1, \qquad 
$$
$$
\lambda_{4}\lambda_{3}^{-1}=1, \qquad 
\lambda_{5}^{-1}\lambda_{4}^{-1}=1, \qquad 
\lambda_{11}^{-1}\lambda_{10} = 1, \qquad 
\lambda_{12}^{-1}\lambda_{14} = 1
$$

Let us simplify these equations. From the last five equations we get:
$$
\rho = \lambda_{3} = \lambda_{4} = \lambda_{5}^{-1} = \lambda_{7}; \qquad 
a = \lambda_{10} = \lambda_{11}; \qquad 
x = \lambda_{12} = \lambda_{14}
$$
Substituting this to the first two equations we get an 
equivalent presentation of our group: 
the generators are $a,x,\beta$; they satisfy two relations: 
\begin{equation} \label{5.3.05.2}
\rho^{-1}a\beta x^{-1}\rho^{-1} =1,\quad \rho^{-1}x^{-1}\rho a  \rho^{-1} =1. 
\end{equation}
Set $b = \rho x^{-1}\rho^{-1}$, i.e. $x^{-1}= \rho^{-1}b\rho$. 
Then the first equation in (\ref{5.3.05.2}) is equivalent to $\beta = ab$. 
Thus the group is generated by $a,b$. The only relation comes from 
the second equation in (\ref{5.3.05.2}). So substituting 
the above expressions for $\beta$ and $x^{-1}$,  we arrive at 
\begin{equation} \label{5.3.05.3}
(ab)^{-1}(ab)^{-1}b (ab) aba (ab)^{-1}=1 \quad \Leftrightarrow \quad bababa = ababab,
\end{equation}
which is the defining relation for the braid group of type $G_2$. 
The proof of Theorem \ref{3.9.05.1} is finished. 

\vskip 3mm
\paragraph{4.4 The action of the braid group of type $G_2$.} This 
group is generated by two elements $a,b$ subject to the single relation 
$ababab= bababa$. 
Let $a$ be the composition of the three
 mutations at the vertices $b_1$, $c_1$, $b_1$, 
and $b$ be the composition of the three mutations at 
the vertices $c_1$, $b_1$, $c_1$. 
Then one checks that they satisfy the above relation. So they are 
generators of the braid group of type $G_2$. 

\vskip 3mm
\paragraph{4.5 $S^3-L$ and the 
discriminant variety for the Coxeter group of type $G_2$.} Let 
 $W_{G_2}$ be the Coxeter group  of type $G_2$. 
It is isomorphic to 
the dihedral group of order $12$. It acts on 
the two dimensional complex vector space $V_2$, the Cartan subalgebra of the 
complex Lie algebra of type $G_2$, equipped with a configuration 
of six one dimensional subspaces, corresponding to 
the kernels of the roots. The group $W_{G_2}$ acts freely on  
the complement $V_2^{\rm reg}$  to the union of this $6$ lines. 
 The quotient 
$V_2^{\rm reg}/W_{G_2}$ is known to be a $K(\pi, 1)$ space for 
$\pi = W_{G_2}$. 
The 
group $\R_+^*$ acts by dilotations on $V_2^{\rm reg}$, commuting with the 
$W_{G_2}$-action. 
Hence $V_2^{\rm reg}/R^*_+W_{G_2}$ is also
 a $K(\pi, 1)$ space for 
$\pi = W_{G_2}$.

\begin{conjecture} \label{3.16.05.1} The space 
$V_2^{\rm reg}/R^*_+W_{G_2}$ is homeomorphic to $S^3-L$. 
\end{conjecture}

Here is an evidence. 
The quotient of $V_2$ by the action of the group $W_{G_2}$ 
is isomorphic to $\C^2$. The $\R^*_+$-action on $V_2$ descends to an
 $\R^*_+$-action on $V_2^{\rm reg}/W_{G_2} = \C^2$ given by 
$t: (z_1, z_2) \lms (t^2z_1, t^6 z_2)$.  Thus the quotient space 
$V_2^{\rm reg}-\{0\}/W_{G_2}\R^*_+$ is a sphere $S^3$, 
given  as the quotient of $\C^2 -\{0.0\}$ 
by the $(2,6)$-weighted $\R^*_+$-action.   

The intersecting each line with the unit sphere in $V_2$ is a circle. 
The action of the group $W_{G_2}$ 
has two orbits on the set of the six lines, and hence on the set of six circles. 
So we get two circles in the quotient. They  are  
the two connected component of the link $L$ in $S^3$.

\section{Appendix 2: Cluster structure of the moduli of  
triples of flags in $PGL_4$}
\vskip 3mm
\paragraph{5.1 The cluster ${\cal X}$-structure for 
the moduli space ${\cal M}_{0, n+3}$.} Recall that  ${\cal M}_{0, n+3}$ 
is the moduli space of configurations of $n+3$ distinct points on ${\Bbb P}^1$. 
Following Section 8 of \cite{FG1}, we provide it with a structure of 
a cluster ${\cal X}$-variety as follows. 

Let us consider an 
 $(n+3)$-gon whose vertices are labeled by a configuration $(x_1, ..., x_{n+3})$ of points 
on ${\Bbb P}^1$, so that the cyclic structure on the points coincides with the one of the vertices induced by the counterclockwise orientation of the $(n+3)$-gon. 
Then given a triangulation $T$ of the
$(n+3)$-gon, we define a rational coordinate system on  ${\cal M}_{0, n+3}$ as follows. 
Recall the cross-ratio of four distinct points $(x_1, x_2, x_3, x_4)$  
on ${\Bbb P}^1$:
$$
r^+(x_1, x_2, x_3, x_4):= \frac{(x_1-x_2)(x_3-x_4)}{(x_1-x_4)(x_2-x_3)}
$$
We assign to every (internal) edge of the triangulation $T$ a rational function 
$X^T_E$ on ${\cal M}_{0, n+3}$ as follows. We define a seed 
$\mathbf I^T = (I^T, I^T_0, \varepsilon_{ij}, d_i)$ as follows: 
$I^T$ is the set of the edges of $T$, $I^T_0$ is empty, and $d_i=1$. 

Let $(x_a, x_b, x_c, x_d)$ 
be the configuration of points assigned to the vertices of the $4$-gon 
formed by the two triangles of the triangulation sharing the edge $E$, so 
that $E = x_bx_d$. Then 
$$
X^T_E(x_1, ... , x_{n+3}):= r^+(x_a, x_b, x_c, x_d)
$$
\begin{figure}[ht]
\centerline{\epsfbox{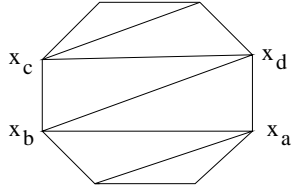}}
\caption{A triangulation of octagon provides a coordinate system  on the moduli space 
${\cal M}_{0,8}$.}
\label{cox5}
\end{figure}

\begin{proposition} \label{2.25.05.4}
The above construction provides ${\cal M}_{0, n+3}$ with a structure of cluster 
${\cal X}$-variety of finite type $A_n$. 
\end{proposition}

{\bf Proof}. Take the snake triangulation shown on Fig \ref{cox5}. 
The corresponding $\varepsilon_{ij}$-function is the one assigned 
to the root system of type $A_n$. 
Let us change a triangulation by flipping a diagonal. Then it is easy to see that 
the resulting transformation of the coordinates is described by the cluster mutation corresponding 
to the flip 
(Section 8 of \cite{FG2}). The proposition is proved.

\vskip 3mm 
\paragraph{5.2 The cluster ${\cal X}$-variety corresponding 
to the moduli space ${\rm Conf}_3({\cal B}_{A_3})$ of triples of flags in $PGL_4$.}
Recall the rational coordinates on the moduli space ${\rm Conf}_3({\cal B}_{A_3})$ 
defined in Section 8.3 of \cite{FG1}. A triple of 
vectors $(a_1, a_2, a_3)$ in a four-dimensional vector space $V_4$ provides a flag 
$(A_0, A_1, A_2):= (a_1, a_1a_2, a_1a_2a_3)$ in ${\Bbb P}(V_4)$. Here  
$A_i$ is the projectivisation of the 
subspace spanned by $a_1$, ..., $a_{i+1}$.  
Consider a a triple of flags in ${\Bbb P}(V_4)$:
$$
(A, B, C) = \Bigl((a_1, a_2, a_3), (b_1, b_2, b_3), (c_1, c_2, c_3)\Bigr) 
$$
Let us choose a volume form $\omega \in {\rm det}V_4^*$. Then for any four vectors 
$a,b,c,d)$ there is a determinant
$$
\Delta(a,b,c,d) := \langle\omega, a\wedge b \wedge c \wedge d\rangle 
$$
We define a rational function $X_1$ on ${\rm Conf}_3({\cal B}_{A_3})$ as follows:
$$
X_1(A, B, C):= -\frac{\Delta(a_1,a_2,a_3, b_1 )\Delta(a_1,b_2, b_3, c_1)
\Delta(a_1, c_1, c_2, a_2)}
{\Delta(a_1,a_2,a_3, c_1  )\Delta(a_1, b_1, b_2, a_2)\Delta(a_1, c_1, c_2, b_1)}
$$
and the functions $X_2$ and $X_3$ are obtained by cyclic shifts:
$$
X_2(A, B, C):= X_1(B, C, A), \quad X_3(A, B, C):= X_1(C, A, B), 
$$

\vskip 3mm 
\paragraph{5.3 An isomorphism of cluster ${\cal X}$-varieties.} Let us define a map 
$$
\Phi: {\cal M}_{0,6} \lra {\rm Conf}_3({\cal B}_{A_3})
$$

Recall that a {\it normal curve} $N \subset {\Bbb P}^n$ 
is a curve of the minimal possible degree $n$ which does not lie in a hyperplane. Any such curve is projectively equivalent to the image of the map $t \mapsto (1,t,t^2,\ldots,t^n-1)$. The group of projective transformations preserving a normal curve  is isomorphic to $PGL_2$, that is to the automorphism group of   ${\Bbb P}^1$. For any $n+3$ generic points in ${\Bbb P}^n$ there exists a unique normal curve containing these points. 

Let $(x_1, y_1, x_2, y_2, x_3, y_3)$ be a configuration of six distinct points on ${\Bbb P}^1$. We identify it with a configuration of points on a normal curve $N\subset {\Bbb P}^3$. Set
$$
\Phi(x_1, y_1, x_2, y_2, x_3, y_3) = (X, Y, Z):= 
\Bigl((x_1, x_1y_1, y_3x_1y_1), (x_2, x_2y_2, y_1x_2y_2), 
(x_3, x_3y_3, y_2x_3y_3)\Bigr) 
$$
\begin{figure}[ht]
\centerline{\epsfbox{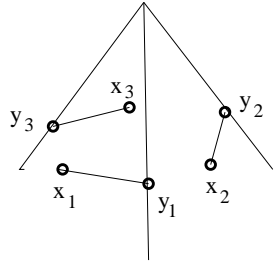}}
\caption{A configuration of three flags in ${\Bbb P}^3$ from a configuration of six points on ${\Bbb P}^1$.}
\label{cox4}
\end{figure}
The inverse map $\Psi$ is defined as follows. Let
\begin{equation} \label{2.25.05.1}
(A, B, C):= (A_0, A_1, A_2), (B_0, B_1, B_2), (C_0, C_1, C_2) 
\end{equation}
be a triple of flags in ${\Bbb P}^3$. So $A_0$ is a point, $A_1$ is a 
line containing this point, and $A_2$ is a plane containing $A_1$. We assign to it the following collection of $6$ points in ${\Bbb P}^3$:
$$
\Psi(A, B, C) :=
(A_0, A_1\cap B_2, B_0, B_1\cap C_2, C_0, C_1\cap A_2)
$$
Taking the unique normal curve passing through these points, we get a configuration of $6$ points on ${\Bbb P}^1$. By the very definition, the compositions  
$\Psi\circ \Phi$ and $\Phi\circ \Psi$ are the identity maps.

\begin{proposition} \label{1.25.05.2}
The map $\Phi$ provides an isomorphism of the cluster ${\cal X}$-varieties corresponding to 
the moduli spaces  ${\cal M}_{0,6}$ and ${\rm Conf}_3({\cal B}_{A_3})$. 
\end{proposition}

\begin{figure}[ht]
\centerline{\epsfbox{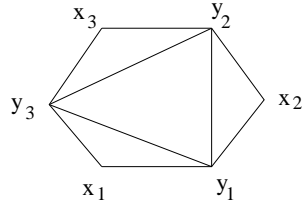}}
\caption{A triangulation of hexagon providing a coordinate system  on 
${\cal M}_{0,6}$.}
\label{cox3}
\end{figure}

{\bf Proof}. Observe that one has 
$$
X_1\Phi(x_1, y_1, x_2, y_2, x_3, y_3) = X_1(X, Y, Z)= 
-\frac{\Delta(x_1, y_1, y_3, x_2)\Delta(x_1, x_2, y_2, x_3)
\Delta(x_1, x_3, y_3, y_1)}{\Delta(x_1, y_1, y_3, x_3)\Delta(x_1, x_2, y_2, y_1)
\Delta(x_1, x_3, y_3, x_2)} 
$$
$$
= \frac{\Delta(x_1, x_2, y_1, y_3)\Delta(x_1, x_2,  x_3, y_2)}
{\Delta(x_1, x_2, y_1, y_2)\Delta(x_1, x_2, y_3, x_3)} = r^+(y_1, y_3, x_3, y_2) 
$$
Using the cyclic shifts, we get
$$
X_2\Phi(x_1, y_1, x_2, y_2, x_3, y_3) = r^+(y_2,y_1,x_1, y_3), \quad 
 X_2\Phi(x_1, y_1, x_2, y_2, x_3, y_3) = r^+(y_3,y_2, x_2, y_1) 
$$
These are the coordinates on ${\cal M}_{0,6}$ assigned to the triangulation of 
the hexagon shown on Fig \ref{cox3}. The proposition is proved. 

{\bf Exercise.} Using the above results, show that the 
cluster ${\cal X}$-variety corresponding to the 
moduli space ${\rm Conf}_3({\cal B}_{B_2})$ 
of triples of flags in $Sp_4$ is of finite type $B_2$. Hint: 
use the triangulations of the hexagon symmetric with respect to the center.

\end{document}